\input amstex



\catcode`\X=12\catcode`\@=11

\def\n@wcount{\alloc@0\count\countdef\insc@unt}
\def\n@wwrite{\alloc@7\write\chardef\sixt@@n}
\def\n@wread{\alloc@6\read\chardef\sixt@@n}
\def\r@s@t{\relax}\def\v@idline{\par}\def\@mputate#1/{#1}
\def\l@c@l#1X{\firstpart.#1}\def\gl@b@l#1X{#1}\def\t@d@l#1X{{}}

\def\crossrefs#1{\ifx\all#1\let\tr@ce=\all\else\def\tr@ce{#1,}\fi
   \n@wwrite\cit@tionsout\openout\cit@tionsout=\jobname.cit 
   \write\cit@tionsout{\tr@ce}\expandafter\setfl@gs\tr@ce,}
\def\setfl@gs#1,{\def\@{#1}\ifx\@\empty\let\next=\relax
   \else\let\next=\setfl@gs\expandafter\xdef
   \csname#1tr@cetrue\endcsname{}\fi\next}
\def\m@ketag#1#2{\expandafter\n@wcount\csname#2tagno\endcsname
     \csname#2tagno\endcsname=0\let\tail=\all\xdef\all{\tail#2,}
   \ifx#1\l@c@l\let\tail=\r@s@t\xdef\r@s@t{\csname#2tagno\endcsname=0\tail}\fi
   \expandafter\gdef\csname#2cite\endcsname##1{\expandafter
     \ifx\csname#2tag##1\endcsname\relax?\else\csname#2tag##1\endcsname\fi
     \expandafter\ifx\csname#2tr@cetrue\endcsname\relax\else
     \write\cit@tionsout{#2tag ##1 cited on page \folio.}\fi}
   \expandafter\gdef\csname#2page\endcsname##1{\expandafter
     \ifx\csname#2page##1\endcsname\relax?\else\csname#2page##1\endcsname\fi
     \expandafter\ifx\csname#2tr@cetrue\endcsname\relax\else
     \write\cit@tionsout{#2tag ##1 cited on page \folio.}\fi}
   \expandafter\gdef\csname#2tag\endcsname##1{\expandafter
      \ifx\csname#2check##1\endcsname\relax
      \expandafter\xdef\csname#2check##1\endcsname{}%
      \else\immediate\write16{Warning: #2tag ##1 used more than once.}\fi
      \multit@g{#1}{#2}##1/X%
      \write\t@gsout{#2tag ##1 assigned number \csname#2tag##1\endcsname\space
      on page \number\count0.}%
   \csname#2tag##1\endcsname}}
\def\multit@g#1#2#3/#4X{\def\t@mp{#4}\ifx\t@mp\empty%
      \global\advance\csname#2tagno\endcsname by 1 
      \expandafter\xdef\csname#2tag#3\endcsname
      {#1\number\csname#2tagno\endcsnameX}%
   \else\expandafter\ifx\csname#2last#3\endcsname\relax
      \expandafter\n@wcount\csname#2last#3\endcsname
      \global\advance\csname#2tagno\endcsname by 1 
      \expandafter\xdef\csname#2tag#3\endcsname
      {#1\number\csname#2tagno\endcsnameX}
      \write\t@gsout{#2tag #3 assigned number \csname#2tag#3\endcsname\space
      on page \number\count0.}\fi
   \global\advance\csname#2last#3\endcsname by 1
   \def\t@mp{\expandafter\xdef\csname#2tag#3/}%
   \expandafter\t@mp\@mputate#4\endcsname
   {\csname#2tag#3\endcsname\lastpart{\csname#2last#3\endcsname}}\fi}
\def\t@gs#1{\def\all{}\m@ketag#1e\m@ketag#1s\m@ketag\t@d@l p
   \m@ketag\gl@b@l r \n@wread\t@gsin
   \openin\t@gsin=\jobname.tgs \re@der \closein\t@gsin
   \n@wwrite\t@gsout\openout\t@gsout=\jobname.tgs }
\outer\def\localtags{\t@gs\l@c@l}
\outer\def\globaltags{\t@gs\gl@b@l}
\outer\def\newlocaltag#1{\m@ketag\l@c@l{#1}}
\outer\def\newglobaltag#1{\m@ketag\gl@b@l{#1}}

\newif\ifpr@ 
\def\m@kecs #1tag #2 assigned number #3 on page #4.%
   {\expandafter\gdef\csname#1tag#2\endcsname{#3}
   \expandafter\gdef\csname#1page#2\endcsname{#4}
   \ifpr@\expandafter\xdef\csname#1check#2\endcsname{}\fi}
\def\re@der{\ifeof\t@gsin\let\next=\relax\else
   \read\t@gsin to\t@gline\ifx\t@gline\v@idline\else
   \expandafter\m@kecs \t@gline\fi\let \next=\re@der\fi\next}
\def\pretags#1{\pr@true\pret@gs#1,,}
\def\pret@gs#1,{\def\@{#1}\ifx\@\empty\let\n@xtfile=\relax
   \else\let\n@xtfile=\pret@gs \openin\t@gsin=#1.tgs \message{#1} \re@der 
   \closein\t@gsin\fi \n@xtfile}

\newcount\sectno\sectno=0\newcount\subsectno\subsectno=0
\newif\ifultr@local \def\ultralocal{\ultr@localtrue}
\def\firstpart{\number\sectno}
\def\lastpart#1{\ifcase#1 \or a\or b\or c\or d\or e\or f\or g\or h\or 
   i\or k\or l\or m\or n\or o\or p\or q\or r\or s\or t\or u\or v\or w\or 
   x\or y\or z \fi}

\def\resetall{\global\advance\sectno by 1\subsectno=0
   \gdef\firstpart{\number\sectno}\r@s@t}
\def\resetsub{\global\advance\subsectno by 1
   \gdef\firstpart{\number\sectno.\number\subsectno}\r@s@t}
\def\newsection#1\par{\resetall\vskip0pt plus.3\vsize\penalty-250
   \vskip0pt plus-.3\vsize\bigskip\bigskip
   \message{#1}\leftline{\bf#1}\nobreak\bigskip}
\def\subsection#1\par{\ifultr@local\resetsub\fi
   \vskip0pt plus.2\vsize\penalty-250\vskip0pt plus-.2\vsize
   \bigskip\smallskip\message{#1}\leftline{\bf#1}\nobreak\medskip}

\def\t@gsoff#1,{\def\@{#1}\ifx\@\empty\let\next=\relax\else\let\next=\t@gsoff
   \def\@@{p}\ifx\@\@@\else
   \expandafter\gdef\csname#1cite\endcsname##1{\zeigen{##1}}
   \expandafter\gdef\csname#1page\endcsname##1{?}
   \expandafter\gdef\csname#1tag\endcsname##1{\zeigen{##1}}\fi\fi\next}
\def\verbatimtags{\ifx\all\relax\else\expandafter\t@gsoff\all,\fi}
\def\zeigen#1{\hbox{$\langle$}#1\hbox{$\rangle$}}

\def\(#1){\edef\dot@g{\ifmmode\ifinner(\hbox{\noexpand\etag{#1}})
   \else\noexpand\eqno(\hbox{\noexpand\etag{#1}})\fi
   \else(\noexpand\ecite{#1})\fi}\dot@g}

\newif\ifbr@ck
\def\eat#1{}
\def\[#1]{\br@cktrue[\br@cket#1'X]}
\def\br@cket#1'#2X{\def\temp{#2}\ifx\temp\empty\let\next\eat
   \else\let\next\br@cket\fi
   \ifbr@ck\br@ckfalse\br@ck@t#1,X\else\br@cktrue#1\fi\next#2X}
\def\br@ck@t#1,#2X{\def\temp{#2}\ifx\temp\empty\let\neext\eat
   \else\let\neext\br@ck@t\def\temp{,}\fi
   \def\teemp{#1}\ifx\teemp\empty\else\rcite{#1}\fi\temp\neext#2X}
\def\resetbr@cket{\gdef\[##1]{[\rtag{##1}]}}
\def\references{\resetbr@cket\newsection References\par}

\newtoks\symb@ls\newtoks\s@mb@ls\newtoks\p@gelist\n@wcount\ftn@mber
    \ftn@mber=1\newif\ifftn@mbers\ftn@mbersfalse\newif\ifbyp@ge\byp@gefalse
\def\defm@rk{\ifftn@mbers\n@mberm@rk\else\symb@lm@rk\fi}
\def\n@mberm@rk{\xdef\m@rk{{\the\ftn@mber}}%
    \global\advance\ftn@mber by 1 }
\def\rot@te#1{\let\temp=#1\global#1=\expandafter\r@t@te\the\temp,X}
\def\r@t@te#1,#2X{{#2#1}\xdef\m@rk{{#1}}}
\def\b@@st#1{{$^{#1}$}}\def\str@p#1{#1}
\def\symb@lm@rk{\ifbyp@ge\rot@te\p@gelist\ifnum\expandafter\str@p\m@rk=1 
    \s@mb@ls=\symb@ls\fi\write\f@nsout{\number\count0}\fi \rot@te\s@mb@ls}
\def\byp@ge{\byp@getrue\n@wwrite\f@nsin\openin\f@nsin=\jobname.fns 
    \n@wcount\currentp@ge\currentp@ge=0\p@gelist={0}
    \re@dfns\closein\f@nsin\rot@te\p@gelist
    \n@wread\f@nsout\openout\f@nsout=\jobname.fns }
\def\m@kelist#1X#2{{#1,#2}}
\def\re@dfns{\ifeof\f@nsin\let\next=\relax\else\read\f@nsin to \f@nline
    \ifx\f@nline\v@idline\else\let\t@mplist=\p@gelist
    \ifnum\currentp@ge=\f@nline
    \global\p@gelist=\expandafter\m@kelist\the\t@mplistX0
    \else\currentp@ge=\f@nline
    \global\p@gelist=\expandafter\m@kelist\the\t@mplistX1\fi\fi
    \let\next=\re@dfns\fi\next}
\def\symbols#1{\symb@ls={#1}\s@mb@ls=\symb@ls} 
\def\bigsymbol{\textstyle}
\symbols{\bigsymbol\ast,\dagger,\ddagger,\sharp,\flat,\natural,\star}
\def\ftnumbers{\ftn@mberstrue} \def\ftsymbols{\ftn@mbersfalse}
\def\paginal{\byp@ge} \def\resetftnumbers{\ftn@mber=1}
\def\ftnote#1{\defm@rk\expandafter\expandafter\expandafter\footnote
    \expandafter\b@@st\m@rk{#1}}

\long\def\jump#1\endjump{}
\def\ssum{\mathop{\lower .1em\hbox{$\textstyle\Sigma$}}\nolimits}

\def\qed{\nobreak\kern 1em \vrule height .5em width .5em depth 0em}
\def\newneq{\hbox{\rlap{\hbox to 1\wd9{\hss$=$\hss}}\raise .1em 
   \hbox to 1\wd9{\hss$\scriptscriptstyle/$\hss}}}
\def\subsetne{\setbox9 = \hbox{$\subset$}\mathrel{\hbox{\rlap
   {\lower .4em \newneq}\raise .13em \hbox{$\subset$}}}}
\def\supsetne{\setbox9 = \hbox{$\subset$}\mathrel{\hbox{\rlap
   {\lower .4em \newneq}\raise .13em \hbox{$\supset$}}}}

\def\vbar{\mathchoice{\vrule height6.3ptdepth-.5ptwidth.8pt\kern-.8pt}
   {\vrule height6.3ptdepth-.5ptwidth.8pt\kern-.8pt}
   {\vrule height4.1ptdepth-.35ptwidth.6pt\kern-.6pt}
   {\vrule height3.1ptdepth-.25ptwidth.5pt\kern-.5pt}}
\def\f@dge{\mathchoice{}{}{\mkern.5mu}{\mkern.8mu}}
\def\b@c#1#2{{\rm \mkern#2mu\vbar\mkern-#2mu#1}}
\def\b@b#1{{\rm I\mkern-3.5mu #1}}
\def\b@a#1#2{{\rm #1\mkern-#2mu\f@dge #1}}
\def\bb#1{{\count4=`#1 \advance\count4by-64 \ifcase\count4\or\b@a A{11.5}\or
   \b@b B\or\b@c C{5}\or\b@b D\or\b@b E\or\b@b F \or\b@c G{5}\or\b@b H\or
   \b@b I\or\b@c J{3}\or\b@b K\or\b@b L \or\b@b M\or\b@b N\or\b@c O{5} \or
   \b@b P\or\b@c Q{5}\or\b@b R\or\b@a S{8}\or\b@a T{10.5}\or\b@c U{5}\or
   \b@a V{12}\or\b@a W{16.5}\or\b@a X{11}\or\b@a Y{11.7}\or\b@a Z{7.5}\fi}}

\catcode`\X=11 \catcode`\@=12


\sectno=1   
\localtags
\newbox\noforkbox \newdimen\forklinewidth
\forklinewidth=0.3pt   
\setbox0\hbox{$\textstyle\bigcup$}
\setbox1\hbox to \wd0{\hfil\vrule width \forklinewidth depth \dp0
                        height \ht0 \hfil}
\wd1=0 cm
\setbox\noforkbox\hbox{\box1\box0\relax}
\def\unionstick{\mathop{\copy\noforkbox}\limits}
\def\nonfork#1#2_#3{#1\unionstick_{\textstyle #3}#2}
\def\nonforkin#1#2_#3^#4{#1\unionstick_{\textstyle #3}^{\textstyle #4}#2}     
%
\setbox0\hbox{$\textstyle\bigcup$}
\setbox1\hbox to \wd0{\hfil{\sl /\/}\hfil}
\setbox2\hbox to \wd0{\hfil\vrule height \ht0 depth \dp0 width
                                \forklinewidth\hfil}
\wd1=0cm
\wd2=0cm
\newbox\doesforkbox
\setbox\doesforkbox\hbox{\box1\box0\relax}
\def\nunionstick{\mathop{\copy\doesforkbox}\limits}

\def\fork#1#2_#3{#1\nunionstick_{\textstyle #3}#2}
\def\forkin#1#2_#3^#4{#1\nunionstick_{\textstyle #3}^{\textstyle #4}#2}     
\def\cite #1{\rm[#1]}
\NoBlackBoxes
\define\mr{\medskip\roster}
\define\sn{\smallskip\noindent}
\define\mn{\medskip\noindent}
\define\bn{\bigskip\noindent}
\define\ub{\underbar}
\define\wilog{\text{without loss of generality}}
\define\ermn{\endroster\medskip\noindent}

\define\dbcu{\dsize\bigcup}
\define \nl{\newline}
\documentstyle {amsppt}
\topmatter
\title {Few non-minimal types and non-structure}\endtitle
\author {Saharon Shelah \thanks{\null\newline
I thank Alice Leonhardt for the beautiful typing \null\newline
Latest Revision - 99/Feb/5\null\newline
Publication 603} \endthanks} \endauthor
\affil{Institute of Mathematics \\
The Hebrew University \\
Jerusalem, Israel
\medskip
Rutgers University \\
Department of Mathematics \\
New Brunswick, NJ USA} \endaffil
\mn
\abstract 
We pay two debts from \cite{Sh:576}. The main demands little
knowledge from \cite{Sh:576}, just quoting a model theoretic consequence of
the weak diamond.  We assume that ${\frak K}$ has amalgamation in $\lambda$,
and that the minimal types are not dense to get many non-isomorphic models
in $\lambda^+$.  For this also pcf considerations are relevant. 

The minor debt was the use in one point of \cite{Sh:576} of $\lambda \ne
\aleph_0$, it is minor as for this case by \cite{Sh:88} we ``usually"
know more.
\endabstract
\endtopmatter
\document  
\def\renewcommand{\newcommand}	       
\edef\cite{\the\catcode`@}%
\catcode`@ = 11
\let\@oldatcatcode = \cite
\chardef\@letter = 11
\chardef\@other = 12
%
%
%
%
\def\@innerdef#1#2{\edef#1{\expandafter\noexpand\csname #2\endcsname}}%
%
%
\@innerdef\@innernewcount{newcount}%
\@innerdef\@innernewdimen{newdimen}%
\@innerdef\@innernewif{newif}%
\@innerdef\@innernewwrite{newwrite}%
%
%
%
\def\@gobble#1{}%
%
%
%
\ifx\inputlineno\@undefined
   \let\@linenumber = \empty 
\else
   \def\@linenumber{\the\inputlineno:\space}%
\fi
%
%
%
\def\@futurenonspacelet#1{\def\cs{#1}%
   \afterassignment\@stepone\let\@nexttoken=
}%
\begingroup 
\def\\{\global\let\@stoken= }%
\\ 
\endgroup
\def\@stepone{\expandafter\futurelet\cs\@steptwo}%
\def\@steptwo{\expandafter\ifx\cs\@stoken\let\@@next=\@stepthree
   \else\let\@@next=\@nexttoken\fi \@@next}%
\def\@stepthree{\afterassignment\@stepone\let\@@next= }%
%
%
%
\def\@getoptionalarg#1{%
   \let\@optionaltemp = #1%
   \let\@optionalnext = \relax
   \@futurenonspacelet\@optionalnext\@bracketcheck
}%
%
%
\def\@bracketcheck{%
   \ifx [\@optionalnext
      \expandafter\@@getoptionalarg
   \else
      \let\@optionalarg = \empty
      \expandafter\@optionaltemp
   \fi
}%
\def\@@getoptionalarg[#1]{%
   \def\@optionalarg{#1}%
   \@optionaltemp
}%
%
%
%
\def\@nnil{\@nil}%
\def\@fornoop#1\@@#2#3{}%
\def\@for#1:=#2\do#3{%
   \edef\@fortmp{#2}%
   \ifx\@fortmp\empty \else
      \expandafter\@forloop#2,\@nil,\@nil\@@#1{#3}%
   \fi
}%
\def\@forloop#1,#2,#3\@@#4#5{\def#4{#1}\ifx #4\@nnil \else
       #5\def#4{#2}\ifx #4\@nnil \else#5\@iforloop #3\@@#4{#5}\fi\fi
}%
\def\@iforloop#1,#2\@@#3#4{\def#3{#1}\ifx #3\@nnil
       \let\@nextwhile=\@fornoop \else
      #4\relax\let\@nextwhile=\@iforloop\fi\@nextwhile#2\@@#3{#4}%
}%
%
%
%
\@innernewif\if@fileexists
\def\@testfileexistence{\@getoptionalarg\@finishtestfileexistence}%
\def\@finishtestfileexistence#1{%
   \begingroup
      \def\extension{#1}%
      \immediate\openin0 =
         \ifx\@optionalarg\empty\jobname\else\@optionalarg\fi
         \ifx\extension\empty \else .#1\fi
         \space
      \ifeof 0
         \global\@fileexistsfalse
      \else
         \global\@fileexiststrue
      \fi
      \immediate\closein0
   \endgroup
}%
%
%
%
%
\def\bibliographystyle#1{%
   \@readauxfile
   \@writeaux{\string\bibstyle{#1}}%
}%
\let\bibstyle = \@gobble
%
%
\let\bblfilebasename = \jobname
\def\bibliography#1{%
   \@readauxfile
   \@writeaux{\string\bibdata{#1}}%
   \@testfileexistence[\bblfilebasename]{bbl}%
   \if@fileexists
      \nobreak
      \@readbblfile
   \fi
}%
\let\bibdata = \@gobble
%
%
\def\nocite#1{%
   \@readauxfile
   \@writeaux{\string\citation{#1}}%
}%
\@innernewif\if@notfirstcitation
%
%
\def\cite{\@getoptionalarg\@cite}%
%
%
\def\@cite#1{%
   \let\@citenotetext = \@optionalarg
   \printcitestart
   \nocite{#1}%
   \@notfirstcitationfalse
   \@for \@citation :=#1\do
   {%
      \expandafter\@onecitation\@citation\@@
   }%
   \ifx\empty\@citenotetext\else
      \printcitenote{\@citenotetext}%
   \fi
   \printcitefinish
}%
\def\@onecitation#1\@@{%
   \if@notfirstcitation
      \printbetweencitations
   \fi
   \expandafter \ifx \csname\@citelabel{#1}\endcsname \relax
      \if@citewarning
         \message{\@linenumber Undefined citation `#1'.}%
      \fi
      \expandafter\gdef\csname\@citelabel{#1}\endcsname{%
\strut
\vadjust{\vskip-\dp\strutbox
\vbox to 0pt{\vss\parindent0cm \leftskip=\hsize 
\advance\leftskip3mm
\advance\hsize 4cm\strut\openup-4pt 
\rightskip 0cm plus 1cm minus 0.5cm ?  #1 ?\strut}}
         {\tt
            \escapechar = -1
            \nobreak\hskip0pt
            \expandafter\string\csname#1\endcsname
            \nobreak\hskip0pt
         }%
      }%
   \fi
   \csname\@citelabel{#1}\endcsname
   \@notfirstcitationtrue
}%
%
%
\def\@citelabel#1{b@#1}%
%
%
\def\@citedef#1#2{\expandafter\gdef\csname\@citelabel{#1}\endcsname{#2}}%
%
%
%
\def\@readbblfile{%
   \ifx\@itemnum\@undefined
      \@innernewcount\@itemnum
   \fi
   \begingroup
      \def\begin##1##2{%
         \setbox0 = \hbox{\biblabelcontents{##2}}%
         \biblabelwidth = \wd0
      }%
      \def\end##1{}
      %
      %
      \@itemnum = 0
      \def\bibitem{\@getoptionalarg\@bibitem}%
      \def\@bibitem{%
         \ifx\@optionalarg\empty
            \expandafter\@numberedbibitem
         \else
            \expandafter\@alphabibitem
         \fi
      }%
      \def\@alphabibitem##1{%
         \expandafter \xdef\csname\@citelabel{##1}\endcsname {\@optionalarg}%
         \ifx\biblabelprecontents\@undefined
            \let\biblabelprecontents = \relax
         \fi
         \ifx\biblabelpostcontents\@undefined
            \let\biblabelpostcontents = \hss
         \fi
         \@finishbibitem{##1}%
      }%
      \def\@numberedbibitem##1{%
         \advance\@itemnum by 1
         \expandafter \xdef\csname\@citelabel{##1}\endcsname{\number\@itemnum}%
         \ifx\biblabelprecontents\@undefined
            \let\biblabelprecontents = \hss
         \fi
         \ifx\biblabelpostcontents\@undefined
            \let\biblabelpostcontents = \relax
         \fi
         \@finishbibitem{##1}%
      }%
      \def\@finishbibitem##1{%
         \biblabelprint{\csname\@citelabel{##1}\endcsname}%
         \@writeaux{\string\@citedef{##1}{\csname\@citelabel{##1}\endcsname}}%
         \ignorespaces
      }%
      %
      %
      \let\em = \bblem
      \let\newblock = \bblnewblock
      \let\sc = \bblsc
      \frenchspacing
      \clubpenalty = 4000 \widowpenalty = 4000
      \tolerance = 10000 \hfuzz = .5pt
      \everypar = {\hangindent = \biblabelwidth
                      \advance\hangindent by \biblabelextraspace}%
      \bblrm
      \parskip = 1.5ex plus .5ex minus .5ex
      \biblabelextraspace = .5em
      \bblhook
      \input \bblfilebasename.bbl
   \endgroup
}%
%
%
\@innernewdimen\biblabelwidth
\@innernewdimen\biblabelextraspace
%
%
%
\def\biblabelprint#1{%
   \noindent
   \hbox to \biblabelwidth{%
      \biblabelprecontents
      \biblabelcontents{#1}%
      \biblabelpostcontents
   }%
   \kern\biblabelextraspace
}%
%
%
%
\def\biblabelcontents#1{{\bblrm [#1]}}%
%
%
\def\bblrm{\rm}%
%
%
\def\bblem{\it}%
%
%
\def\bblsc{\ifx\@scfont\@undefined
              \font\@scfont = cmcsc10
           \fi
           \@scfont
}%
%
%
\def\bblnewblock{\hskip .11em plus .33em minus .07em }%
%
%
\let\bblhook = \empty
%
%
%
\def\printcitestart{[}
\def\printcitefinish{]}
\def\printbetweencitations{, }
\def\printcitenote#1{, #1}
%
%
%
\let\citation = \@gobble
%
%
%
\@innernewcount\@numparams
%
%
\def\newcommand#1{%
   \def\@commandname{#1}%
   \@getoptionalarg\@continuenewcommand
}%
%
%
\def\@continuenewcommand{%
   \@numparams = \ifx\@optionalarg\empty 0\else\@optionalarg \fi \relax
   \@newcommand
}%
%
%
\def\@newcommand#1{%
   \def\@startdef{\expandafter\edef\@commandname}%
   \ifnum\@numparams=0
      \let\@paramdef = \empty
   \else
      \ifnum\@numparams>9
         \errmessage{\the\@numparams\space is too many parameters}%
      \else
         \ifnum\@numparams<0
            \errmessage{\the\@numparams\space is too few parameters}%
         \else
            \edef\@paramdef{%
               \ifcase\@numparams
                  \empty  No arguments.
               \or ####1%
               \or ####1####2%
               \or ####1####2####3%
               \or ####1####2####3####4%
               \or ####1####2####3####4####5%
               \or ####1####2####3####4####5####6%
               \or ####1####2####3####4####5####6####7%
               \or ####1####2####3####4####5####6####7####8%
               \or ####1####2####3####4####5####6####7####8####9%
               \fi
            }%
         \fi
      \fi
   \fi
   \expandafter\@startdef\@paramdef{#1}%
}%
%
%
%
%
\def\@readauxfile{%
   \if@auxfiledone \else 
      \global\@auxfiledonetrue
      \@testfileexistence{aux}%
      \if@fileexists
         \begingroup
            \endlinechar = -1
            \catcode`@ = 11
            \input \jobname.aux
         \endgroup
      \else
         \message{\@undefinedmessage}%
         \global\@citewarningfalse
      \fi
      \immediate\openout\@auxfile = \jobname.aux
   \fi
}%
%
%
\newif\if@auxfiledone
\ifx\noauxfile\@undefined \else \@auxfiledonetrue\fi
%
%
%
%
\@innernewwrite\@auxfile
\def\@writeaux#1{\ifx\noauxfile\@undefined \write\@auxfile{#1}\fi}%
%
%
%
\ifx\@undefinedmessage\@undefined
   \def\@undefinedmessage{No .aux file; I won't give you warnings about
                          undefined citations.}%
\fi
%
%
\@innernewif\if@citewarning
\ifx\noauxfile\@undefined \@citewarningtrue\fi
%
%
%
\catcode`@ = \@oldatcatcode
\def\widestnumber#1#2{}
\def\rm{\fam0 \tenrm}
\def\fakesubhead#1\endsubhead{\bigskip\noindent{\bf#1}\par}
\newpage

In \cite{Sh:576} there was one point where we used as assumption $I(\lambda
^{+3},K) = 0$.  This was fine for the purpose there, but is unsuitable in
our present framework: we want to analyze what occurs in higher cardinals, so 
our main aim here is to eliminate its use and add to our knowledge on
non-structure.

The point was ``the minimal triples in $K^3_\lambda$ are dense" 
(\cite{Sh:576},3.17).  For this we assume we have a counterexample, and try
to build many nonisomorphic models; 
hence we get cases of amalgamation which are
necessarily unique.  So we try to build many models in $\lambda^+$ by omitting
``types" over models of size $\lambda$, in a specific way where unique
amalgamation holds.  If this argument fails, we prove $\bold C^1_{{\frak K},
\lambda}$ has weak $\lambda^+$-coding and by it get $2^{\lambda^{++}}$ 
non-isomorphic models except when the weak diamond ideal on $\lambda^+$ 
is $\lambda^{++}$-saturated and use pcf to get the full result.  We work also
to get large IE (many models no one $\le_{\frak K}$-embedding to another).

There was also another point left in \cite[4.2]{Sh:576}, for the case
$\lambda = \aleph_0$ only this is filled in the end.

A third point we needed to prove in \cite{Sh:576} was 
\cite[Th.6.5]{Sh:576}, see notation there: $K^{2,uq}_\lambda = 
\emptyset \Rightarrow I(\lambda^{++},K) = 
2^{\lambda^{++}}$ (assuming $2^n < 2^{\lambda^+} < 2^{\lambda^{++}}$,
etc. we just prove $I(\lambda^{++},K) \ge \mu_{\text{wd}}
(\lambda^+))$.  This holds as the context of \cite{Sh:576} is included
in our context ($\lambda$-good$^-$).
\bn
We assume here some knowledge of \cite[\S2]{Sh:576}.
\demo{\stag{2.0} Context}
\roster
\item "{$(a)$}"  ${\frak K}$ abstract elementary class $LS({\frak K}) \le
\lambda$
\sn
\item "{$(b)$}"  ${\frak K}$ has amalgamation in $\lambda$.
\endroster
\enddemo
\bigskip

\definition{\stag{2.1} Definition}  1) For $x \in \{a,d\}$ we say
$UQ^x_\lambda(M_0,M_1,M_2,M_3)$ if:
\medskip
\roster
\item "{$(a)$}"  $M_\ell \in K_\lambda$ for $\ell \le 3$
\smallskip
\noindent
\item "{$(b)$}"  $M_0 \le_{\frak K} M_\ell$ for $\ell =1,2$
\smallskip
\noindent
\item "{$(c)$}"  if for $i \in \{1,2\}$ we have $M^i_\ell \in K_\lambda$, for
$\ell < 4$ and $M^i_0 \le_{\frak K} M^i_\ell \le_{\frak K} M^i_3$ for 
$i=1,2,\ell =1,2$ and $[x=d \Rightarrow M^i_1 \cap M^i_2 = M^i_0]$ and 
$f^i_\ell$ is an isomorphism from $M_\ell$ onto $M^i_\ell$ for $\ell < 3$ and 
$f^i_0 \subseteq f^i_1,f^i_0 \subseteq f^i_2$ \underbar{then} there are 
$M'_3,f_3$ such that $M^2_3 \le_{\frak K} M'_3$ and $f_3$ is a 
$\le_{\frak K}$-embedding of $M^1_3$ into $M'_3$ extending \newline 
$(f^2_1 \circ (f^1_1)^{-1}) \cup (f^2_2 \circ (f^1_2)^{-1})$ i.e. $f_3 \circ
f^1_1 = f^2_1 \and f_3 \circ f^1_2 = f^2_2$
\sn
\item "{$(d)$}"  $M_0 \le_{\frak K} M_\ell \le_{\frak K} M_3 \in K_\lambda$ 
for $\ell =1,2$
\sn
\item "{$(e)$}"  $x=d \Rightarrow M_1 \cap M_2 = M_0$.
\ermn
2) We say $UQ^x_\lambda(M_0,M_1,M_2)$ if $UQ^x_\lambda(M_0,M_1,M_2,M_3)$ for
some $M_3$. \nl
3) If we omit $x$, we mean $x=a$. \nl
4) $K^{3,*_m}_\lambda$ is the family of triple $(M,N,a) \in K^3_\lambda$ such
that there is no minimal triple above it. \nl
5) $K^{2,*}_\lambda$ is the family $\{(M,N):$ for some $a,(M,N,a) \in
K^{3,*}_\lambda$. \nl
6) ${\Cal S}_*(M) = \{p \in {\Cal S}(M)$: for some $(M,N,a) \in K^{3,*}
_\lambda$ we have $p = \text{ tp}(a,M,N)\}$.
\enddefinition
\bigskip

\proclaim{\stag{2.2} Claim}  1) Symmetry: assuming $x \in \{a,d\}$ we have
$UQ^x_\lambda(M_0,M_1,M_2,M_3) \Rightarrow UQ^x_\lambda(M_0,M_2,M_1,M_3)$;
we can also omit $M_3$. \newline
2) $UQ^a_\lambda(M_0,M_1,M_2) \Rightarrow UQ^d_\lambda(M_0,M_1,M_2)$ recalling
$M_0$ is an amalgamation base (in ${\frak K}_\lambda$). \newline
3) $UQ^a_\lambda(M_0,M_1,M_2,M_3)$ \underbar{iff} clauses (a), (b),
(d), (e) of Definition \scite{2.1}(1),(2) holds and also $(c)^-$, i.e. 
clause (c) restricted to the case $M^1_\ell = M_\ell$ for $\ell \le 3$. \nl
4) If $UQ^x_\lambda(M_0,M_1,M_2,M_3),M_3 \le_{\frak K} M'_3 \in K_\lambda$ 
then $UQ^x_\lambda(M_0,M_1,M_2,M'_3)$; and also the inverse: if
$UQ^x_\lambda(M_0,M_1,M_2,M'_3)$ and $M_1 \cup M_2 \subseteq M_3
\le_{\frak K} M'_3$ then \nl
$UQ^x_\lambda(M_0,M_1,M_2,M_3)$. \nl
5) Assume $(M,N,a) \in K^3_\lambda$ and no triple above it is minimal then
$\neg UQ(M,N,N)$.
\endproclaim
\bigskip

\demo{Proof}  1),2)  Trivial. \newline
3) Chasing arrows, we should prove clause (c) of Definition \scite{2.1}(1).  
Assume we are given 
$\langle M^1_\ell:\ell < 4 \rangle,\langle M^2_\ell:\ell < 4 \rangle,
\langle f^i_\ell:\ell < 3 \rangle$ as there for $i=1,2$.  First for $i=1,2$
apply clause $(c)^-$ to $\langle M^i_\ell:\ell < 4 \rangle,\langle f^i_\ell:
\ell < 3 \rangle$.  So there are $N^i_3,f^i_3$ such that: $M^i_3 \le_{\frak K}
N^i_3 \in K_\lambda$, and $f^i_3$ a $\le_{\frak K}$-embedding of $M_3$ into
$N^i_3$ extending $f^i_1 \cup f^i_2$.  As ${\frak K}$ has amalgamation in
$\lambda$ (by \scite{2.0}(b)) there are $N \in K_\lambda$ and 
$\le_{\frak K}$-embeddings
$g^i:N^i \rightarrow N$ such that $g^1 \circ f^1_3 = g^2 \circ f^2_3$, so
we are done.  \newline
4) Again by the amalgamation i.e. \scite{2.0}(b). \nl
5) By \cite[2.6(1)]{Sh:576}.   \hfill$\square_{\scite{2.2}}$
\enddemo
\bigskip

\proclaim{\stag{2.3} Claim}  1) transitivity: If 
$UQ_\lambda(M_\ell,N_\ell,M_{\ell+1},N_{\ell+1})$ for $\ell =0,1$ 
\ub{then} \newline
$UQ_\lambda(M_0,N_0,M_2,N_2)$. \newline
2) If $\theta = \text{ cf}(\theta) < \lambda^+$, and $\langle M_i:i \le \theta
\rangle$ is $\le_{\frak K}$-increasing continuous and $\langle N_i:i \le
\theta \rangle$ is $\le_{\frak K}$-increasing and $UQ_\lambda(M_i,N_i,
M_{i+1},N_{i+1})$ for each $i < \theta$ \underbar{then} \newline
$UQ_\lambda(M_0,N_0,M_\theta,N_\theta)$. \newline
3) Assume:
\medskip
\roster
\item "{$(a)$}"  $\alpha,\beta < \lambda^+$
\sn
\item "{$(b)$}"  $M_{i,j} \in K_\lambda$ for $i \le \alpha,j \le \alpha$
\sn
\item "{$(c)$}"  $i_1 \le i_2 \le \alpha \and j_1 \le j_2 \le \beta
\Rightarrow M_{i_1,j_1} \le_{\frak K} M_{i_2,j_2}$
\sn
\item "{$(d)$}"  $\langle M_{i,j}:i \le \alpha \rangle$ is $\le_{\frak K}$-
increasing continuous for each $j \le \beta$
\sn
\item "{$(e)$}"  $\langle M_{i,j}:j \le \beta \rangle$ is 
$\le_{\frak K}$-increasing continuous for each $i \le \alpha$
\sn
\item "{$(f)$}"  $UQ_\lambda(M_{i,j},M_{i+1,j},M_{i,j+1},M_{i+1,j+1})$.
\endroster
\medskip

\noindent
\underbar{Then} 
$UQ_\lambda(M_{0,0},M_{\alpha,0},M_{0,\beta},M_{\alpha,\beta})$. \newline
4) If $UQ^x_\lambda(M_0,M_1,M_2)$ and $M_0 \le_{\frak K} M'_1 \le_{\frak K}
M_1$ and $M_0 \le_{\frak K} M'_2 \le_{\frak K} M_2$ \underbar{then}
$UQ_\lambda(M_0,M'_1,M'_2)$.
\endproclaim
\bigskip

\demo{Proof}  Chasing arrows (note: $UQ = UQ^a$ is easier than $UQ^d$, for
$UQ^d$ the parallel claim is not clear at this point, e.g. seemingly
transitivity fails). \hfill$\square_{\scite{2.3}}$
\enddemo
\bigskip

\definition{\stag{2.4} Definition}  1) ${\Cal T}_\lambda[{\frak K}] = 
\{(M,\Gamma):M \in K_\lambda,\Gamma \subseteq K_\lambda,|\Gamma| \le 
\lambda \text{ and}$

$\qquad \qquad \qquad \qquad \qquad \qquad \qquad \qquad
N \in \Gamma \Rightarrow M \le_{\frak K} N\}$. \newline
2) We define relations on ${\Cal T}_\lambda[{\frak K}]$:
\medskip
\roster
\item "{$(\alpha)$}"  $(M_1,\Gamma_1) \le_h (M_2,\Gamma_2)$ \underbar{iff}
{\roster
\itemitem{ (a) }  $M_1 \le_{\frak K} M_2$
\sn
\itemitem{ (b) }  $h$ is a partial function from $\Gamma_2$ onto $\Gamma_1$
\sn
\itemitem{ (c) }  if $h(N_2) = N_1$ then $UQ_\lambda(M_1,N_1,M_2,N_2)$.
\endroster}
\item "{$(\beta)$}"  $(M_1,\Gamma_1) \le^h (M_2,\Gamma_2)$ \underbar{iff}
$h$ is a one to one function from $\Gamma_1$ to $\Gamma_2$ such that \newline
$(M_1,\Gamma_2) \le_{h^{-1}} (M_2,\Gamma_2)$
\mn
\item "{$(\gamma)$}"  $(M_1,\Gamma_1) \le^{h_2}_{h_1} (M_2,\Gamma_2)$
\underbar{iff} $(M_1,\Gamma_1) \le_{h_1} (M_2,\Gamma_2),(M_1,\Gamma_2) 
\le^{h_2} (M_2,\Gamma_2)$ and $h_1 \circ h_2 = \text{ id}_{\Gamma_1}$.
\endroster
\medskip

\noindent
3) We write $<_h$ (or $<^h$ or $<^{h_2}_{h_1}$) if in addition 
$M_1 \ne M_2$ moreover \footnote{note that: if there is no minimal triple in
$K^3_\lambda$ then every $a \in M_2 \backslash M_2$ is as required so the
moreover is not necessary} $(M_1,M_2) \in K^{2,*}_\lambda$. \newline
4) $(M_1,\Gamma_1) \le (M_2,\Gamma_2)$ means:
\medskip
\roster
\item "{$(a)$}"  $\langle N \backslash M_\ell:N \in \Gamma_\ell \rangle$ are
pairwise disjoint non-empty (and $(M_\ell,\Gamma_\ell) \in
{\Cal T}_\lambda[{\frak K}]$) for $\ell =1,2$
\sn
\item "{$(b)$}"  for every $N_1 \in \Gamma_1$ there is exactly one
$N_2 \in \Gamma_2$ such that $(N_1 \backslash M_1) \cap (N_2 \backslash
M_2) \ne \emptyset$ (so $(M_1,\Gamma_1) \le (M_1,\Gamma_1)$ just means
clause (a) and we call $(M_1,\Gamma_1)$ disjoint)
\sn
\item "{$(c)$}"  $(M_1,\Gamma_1) \le^h (M_2,\Gamma_2)$ when 
$h:\Gamma_1 \rightarrow \Gamma_2$ is defined by: $h(N_1) = N_2 \Leftrightarrow
(N_1 \backslash M_1) \cap (N_2 \backslash M_2) \ne \emptyset$.
\endroster
\mn
5) Let ${\Cal T}^{dis}_\lambda[{\frak K}] = \{(M,\Gamma) \in {\Cal T}_\lambda
[{\frak K}]:(M,\Gamma) \le (M,\Gamma)\}$.
\enddefinition
\bigskip

\remark{Remark}  1) Of course, except for bookkeeping only the isomorphic
type of the $N \in \Gamma$ over $M$ matters. \newline
2) In \scite{2.5}(2) below, if the continuity is not demanded, the result may
be that $h_{i,\delta}$ is not onto $\Gamma_i$.  So we may use
$\le^{h_2}_{h_1}$. \nl
3) Note that we can use \scite{2.5}(3) and not \scite{2.5}(2).
\endremark
\bigskip

\proclaim{\stag{2.5} Claim}  1) If $(M_1,\Gamma_1) \le_{h_1} (M_2,\Gamma_2)
\le_{h_2} (M_3,\Gamma_3)$ then \newline
$(M_1,\Gamma_1) \le_{h_1 \circ h_2}
(M_3,\Gamma_3)$. \newline
2) If $\delta < \lambda^+$ is a limit ordinal, for $i < j < \delta,
(M_i,\Gamma_i) \le^{h_{i,j}} (M_j,\Gamma_j)$ and the $h_{i,j}$ commute, and
for limit $\alpha < \delta,(\Gamma_\alpha,\langle h_{i,\alpha}:i < \alpha
\rangle)$ is the (direct) limit of $\langle \Gamma_i,h_{i,j}:i < j < 
\alpha \rangle$
\underbar{then} we can find $(M_\delta,\Gamma_\delta)$ and $\langle 
h_{i,\delta}:i < \delta \rangle$ such that $(M_i,\Gamma_i) \le_{h_{i,\delta}} 
(M_\delta,\Gamma_\delta)$ and $M_\delta = \dsize \bigcup_{i < \delta} M_i$ and
$(\Gamma_\delta,\langle h_{i,\delta}:i < \delta \rangle)$ is the limit of
$\langle \Gamma_i,h_{i,j}:i < j < \delta \rangle$. \newline
3) $({\Cal T}_\lambda[{\frak K}],\le)$ is a partial order in which any
increasing sequence of length $< \lambda^+$ has a lub ($=$ the natural limit).
(Actually, only ${\Cal T}^{dis}_\lambda[{\frak K}]$ matters.)
\endproclaim
\bigskip

\demo{Proof}  1) The main point is the use of \scite{2.3}(1). \newline
2) In addition to chasing arrows we use \scite{2.3}(2) (or \scite{2.3}(3) 
for $\beta=1$). \newline
3) Reformulation of 1), 2).  \hfill$\square_{\scite{2.5}}$
\enddemo
\bigskip

\proclaim{\stag{2.6} Claim}  Assume $2^\lambda < 2^{\lambda^+}$ or at least
the definitional \footnote{we could make the relation on $(M_{\eta \char 94
\langle 0 \rangle},\Gamma_{\eta \char 94 \langle 0 \rangle}),(M_{\eta \char 94
\langle 1 \rangle},\Gamma_{\eta \char 94 \langle 1 \rangle})$ symmetric} 
weak diamond;i.e. DfWD$^+(\lambda^+)$, see \cite[1.7]{Sh:576}.  
If $(*)_\lambda$ or at least $(*)'_\lambda$ below holds (hence above 
some triple from $K^3_\lambda$ there is no minimal one), 
\underbar{then} $I(\lambda^+,K) \ge 2^{\lambda^+}$ where
\medskip
\roster
\item "{$(*)_\lambda$}"  for every $(M,\Gamma) \in 
{\Cal T}_\lambda[{\frak K}]$ for some $h,M',\Gamma'$ we have \newline
$(M,\Gamma) <^h (M',\Gamma')$ (without loss of generality $h$ is onto
$\Gamma'$) or just
\smallskip
\noindent
\item "{$(*)'_\lambda$}"  for some $(M_0,\Gamma_0) \in 
{\Cal T}_\lambda[{\frak K}]$, if $(M_0,\Gamma_0) \le^{h_0} (M,\Gamma)$ 
then for some $h,M',\Gamma'$ we have $(M,\Gamma) <^h (M',\Gamma')$.
\endroster
\endproclaim
\bigskip

\demo{Proof}  We choose by induction on $\alpha < \lambda,\langle
(M_\eta,\Gamma_\eta,\Gamma^+_\eta):\eta \in {}^\alpha 2 \rangle$ such that:
\medskip
\roster
\item "{$(a)$}"  $M_\eta \in K_\lambda$ has universe $\gamma_\eta < 
\lambda^+$
\sn    
\item "{$(b)$}"  $(M_\eta,\Gamma_\eta) \in {\Cal T}_\lambda[{\frak K}]$ 
\sn
\item "{$(c)$}"  $(M_\eta,\Gamma_\eta)$ is disjoint (see \scite{2.4}(4)(b))
and $N \in \Gamma_\eta
\rightarrow (N \backslash M_\eta) \cap \lambda^+ = \emptyset$
\sn
\item "{$(d)$}"  $\nu \triangleleft \eta \Rightarrow (M_\nu,\Gamma_\nu)
< (M_\eta,\Gamma_\eta)$
\sn
\item "{$(e)$}"  $(M_\eta,\Gamma^+_\eta) \in {\Cal T}_\lambda[{\frak K}]$ is
disjoint (see \scite{2.4}(4)(b)) and \newline
$N \in \Gamma^+_\eta \Rightarrow (N \backslash M_\eta) \cap \lambda^+ = 
\emptyset$
\sn
\item "{$(f)$}"  $\Gamma_\eta \subseteq \Gamma^+_\eta$
\sn
\item "{$(g)$}"  $(M_\eta,\Gamma^+_\eta) \le (M_{\eta \char 94 \langle 0
\rangle},\Gamma_{\eta \char 94 \langle 0 \rangle})$
\sn
\item "{$(h)$}"  for some $N \in \Gamma^+_\eta$ we have $N \cong_{M_\eta}
M_{\eta \char 94 <1>}$
\sn
\item "{$(i)$}"  $(M_\eta,M_{\eta \char 94 \langle 1 \rangle}) \in K^{2,*}
_\lambda$, that is for some $a \in M_{\eta \char 94 \langle 1 \rangle}$ we
have $(M_\eta,M_{\eta \char 94 \langle 1 \rangle},a) \in K^3_\lambda$ and
above it there is no minimal triple from $K^3_\lambda$.
\endroster
\medskip

There is no problem to carry the induction with $\Gamma^+_\eta \, (\eta \in
{}^\alpha 2)$ chosen in the $(\alpha +1)$-th step.  For $\alpha = 0$ let
$(M_{<>},\Gamma_{<>})$ be the $(M_0,\Gamma_0)$ from $(*)'_\lambda$ 
except that we rename the elements to make the relevant parts of clauses (a),
(c) true.  For $\alpha$ limit use \scite{2.5}(3) (part on lub).  
For $\alpha = \beta+1,\eta \in {}^\beta 2$, by $(*)'_\lambda$ we can find 
$(M_{\eta \char 94 \langle 1 \rangle},\Gamma_{\eta \char 94 \langle 1 
\rangle})$ such that $(M_{\eta \char 94 \langle 1 \rangle},
\Gamma_{\eta \char 94 \langle 1 \rangle}) \in {\Cal T}_\lambda[{\frak K}],
(M_\eta,\Gamma_\eta) < (M_{\eta \char 94 \langle 1 \rangle},
\Gamma_{\eta \char 94 \langle 1 \rangle})$.  
So by the Definition \scite{2.4}(3) of $<$, for some 
$a \in M_{\eta \char 94 \langle 1 \rangle} \backslash
M_\eta$, no triple in $K^3_\lambda$ above $(M_\eta,M_{\eta \char 94 \langle
1 \rangle},a)$ is minimal.  By renaming without loss 
of generality the universe of $M_{\eta \char 94 \langle 1 \rangle}$ is 
some $\gamma_{\eta \char 94 \langle 1 \rangle}(\in (\gamma_\eta,\lambda^+))$
and clause (c) holds.
Let $N_\eta$ be isomorphic to $M_{\eta \char 94 \langle 1 \rangle}$ over 
$M_\eta$ with universe disjoint to
$\lambda^+ \cup \dsize \bigcup_{N \in \Gamma_\eta} |N| \backslash 
\gamma_\eta$ and let $\Gamma^+_\eta = \Gamma_\eta \cup \{N_\eta\}$, so
$(M_\eta,\Gamma^+_\eta) \in {\Cal T}_\lambda[{\frak K}]$ is disjoint, now
apply to it $(*)'_\lambda$ to get $(M_{\eta \char 94 \langle 0 \rangle},
\Gamma_{\eta \char 94 \langle 0 \rangle})$.  Why does clause (h) hold?  By
the choice of $N_\eta$.  So $M_\eta,\Gamma_\eta,\Gamma^+_\eta \, (\eta \in 
{}^{\lambda^+} 2)$ are defined. \nl
Note: if $\eta \char 94 \langle 0 \rangle \triangleleft \nu \in
{}^{\lambda^+ >}2$, \ub{then} $M_{\eta \char 94 <1>}$ is not 
$\le_{\frak K}$-embeddable into $M_\nu$ over $M_\eta$ (by clause (g) + (h)
because by \scite{2.2}(5) and clause (i) we have $\neg UQ(M_\eta,
M_{\eta \char 94 \langle 1 \rangle},N_\eta))$.
By \cite[1.4]{Sh:576} we get the desired conclusion (really also on IE).
\hfill$\square_{\scite{2.6}}$
\enddemo
\bigskip

\proclaim{\stag{2.7} Claim}  An equivalent condition for $(*)_\lambda$ or 
just $(*)'_\lambda$ of \scite{2.6} is (respectively):
\medskip
\roster
\item "{$(**)_\lambda$}"  for every $M \le_{\frak K} N$ from $K_\lambda$ 
for some $M',M <_{\frak K} M' \in K_\lambda$ and \newline
$UQ_\lambda(M,M',N) \and (M,M') \in K^{2,*}_\lambda$
\endroster
\noindent
or just
\roster
\item "{$(**)'_\lambda$}"  for some 
$M_0 \in K_\lambda$ if $M_0 \le_{\frak K} M \le_{\frak K} N \in K_\lambda$ 
then for some $M'$, \newline
$M <_{\frak K} M' \in K_\lambda$ and $UQ_\lambda(M,M',N) \and (M,M') \in
K^{2,*}_\lambda$.
\endroster
\endproclaim
\bigskip

\demo{Proof}  For any $(M,\Gamma) \in {\Cal T}_\lambda[{\frak K}]$, by
``${\frak K}$ has amalgamation in $\lambda$" (and properties of abstract
elementary classes) there are $N^*,\langle f_N:N \in \Gamma \rangle$ such
that:
\medskip
\roster
\item "{$(a)$}"   $M \le_{\frak K} N^* \in K_\lambda$
\sn
\item "{$(b)$}"  for $N \in \Gamma,f_N$ is a $\le_{\frak K}$-embedding of
$N$ into $N^*$ over $M$.
\endroster
\medskip

This shows $(**)_\lambda \Rightarrow (*)_\lambda$ and also $(**)'_\lambda 
\Rightarrow (*)'_\lambda$.
\newline
The other direction are by applying $(*)_\lambda$ (or $(*)'_\lambda$) 
to $(M,\{N\})$. \hfill$\square_{\scite{2.7}}$
\enddemo
\bigskip

\proclaim{\stag{2.8} Claim}  Assume
\medskip
\roster
\item "{$(a)$}"  $(**)_\lambda$ of \scite{2.7} fails
\sn
\item "{$(b)$}"  $M \in K_\lambda \Rightarrow |{\Cal S}_*(M)| > \lambda^+$
(follows from ``above $(M,N,a) \in K^3_\lambda$ there \newline
is no minimal triple" $+ 2^\lambda > \lambda^+$)
\sn
\item "{$(c)$}"  $K$ is categorical in $\lambda$
\sn
\item "{$(d)$}"  $K$ is categorical in $\lambda^+$
\sn
\item "{$(e)$}"  $2^\lambda < 2^{\lambda^+} < 2^{\lambda^{++}}$.
\ermn
\ub{Then} $I(\lambda^{++},K) = 2^{\lambda^{++}}$ except possibly when 
\mr
\item "{$(*)$}"  WDmId$(\lambda^+)$ is $\lambda^{++}$-saturated 
(normal ideal on $\lambda^+$).
\endroster
\endproclaim
\bigskip

\demo{Proof}  We rely on \cite[\S3]{Sh:576}.  By (a) (and (c)) for every 
$M \in K_\lambda$ for some $N = N_{[M]}$ we have $M <_{\frak K} N \in 
K_\lambda$ and \wilog \, even $(M,N) \in K^{2,*}_\lambda$ 
and $(M,M') \in K^{2,*}_\lambda \Rightarrow \neg UQ(M,M',N)$, 
of course, using $N'$ where $N \le_{\frak K} N' \in K_\lambda$ will do too
as ${\frak K}$ has amalgamation in $\lambda$.

Let $\bold C = \bold C^1_{\lambda^+}[{\frak K}]$, see 
\cite[Definition 3.3]{Sh:576}.  
So there is $\langle M^*_i:i < \lambda^+ \rangle$ such that $(M^*_i,
M^*_{i+1}) \cong (M^*_i,N_{[M^*_i]})$.

Now for $\langle M_i:i < \lambda^+ \rangle \in {\bold Seq}_{\lambda^+}
[\bold C]$, by clause (b) there is $p^* \in {\Cal S}_*(M_0)$ not realized in
$\dbcu_{i < \lambda^+} M_i$ hence $M_i <_K M' \in K_\lambda \and$ (some
$a \in M'$ realizes $p$) $\Rightarrow (M,M') \in K^{2,*}_\lambda$;
for a club of $i$'s, for some $j \in (i,\lambda^+)$ we have, 
by \scite{2.3}(4) as $K$ is categorical in $\lambda^+$:
\mr
\item "{$(\alpha)$}"  for some $\le_{\frak K}$-embedding $f_i$ of
$N_{[M^*_i]}$ into $M_j,f_i$ maps $M^*_i$ onto $M_i$ so by assumption (a)
\sn
\item "{$(\beta)$}"  $M_i <_{\frak K} M' \in K_\lambda \and$ some $a \in M'$ 
realizes $p^* \Rightarrow \neg UQ_\lambda(M_i,M_j,M')$.
\ermn

\relax From this it follows that for some
$\lambda^+$-amalgamation function $F,F$ has the weak $\lambda^+$-coding 
property for $\bold C$ (see \cite[3.8(2)]{Sh:576} hence the weaker version, 
too.  Hence by \cite[3.12]{Sh:576}, if WDmId$(\lambda^+)$ is 
not $\lambda^{++}$-saturated then $I(\lambda^{++},K) = 2^{\lambda^{++}}$.
\newline
${}$ \hfill$\square_{\scite{2.8}}$
\enddemo
\bigskip

\remark{\stag{2.9} Remark}  1) We can get more abstract results. \newline
2) Note $\neg(*)_\lambda$ of \scite{2.8} is a ``light" assumption, in fact,
e.g. its negation has high consistency strength.
\endremark
\bigskip

\centerline {$* \qquad * \qquad *$}
\bigskip

The following will provide us a useful division into cases (it is from pcf
theory; on $\mu_{\text{wd}}(\lambda)$ see \cite[1.1]{Sh:576}), we can replace 
$\lambda^+$ by regular $\lambda$ such that $2^\theta = 
2^{< \lambda} < 2^\lambda$ for some $\theta$).
\bigskip

\demo{\stag{2.10} Fact}   Assume $2^\lambda < 2^{\lambda^+}$.

Then one of the following cases occurs:
\medskip
\roster
\item "{$(A)_\lambda$}"  we can find $\mu$ such that letting $\chi^* =
2^{\lambda^+}$
{\roster
\itemitem{ $(\alpha)$ }  $\lambda^+ < \mu \le 2^\lambda$ and cf$(\mu) 
= \lambda^+$
\sn
\itemitem{ $(\beta)$ }   $pp(\mu) = \chi^*$, moreover $pp^+(\mu) =
(\chi^*)^+$ and $\chi^* > 2^\lambda$
\sn
\itemitem{ $(\gamma)$ }  $(\forall \mu')(\text{cf}(\mu') \le \lambda^+ < \mu'
< \mu \rightarrow pp(\mu') < \mu)$ hence \newline
cf$(\mu') \le \lambda^+ < \mu' < \mu \Rightarrow pp_{\lambda^+}(\mu') < \mu$
\sn
\itemitem{ $(\delta)$ }  for every regular cardinal $\chi$ in the interval
$(\mu,\chi^*]$ there is an increasing sequence $\langle \lambda_i:
i < \lambda^+ \rangle$ of regular cardinals $> \lambda^+$ with limit $\mu$ 
such that $\chi = \text{ tcf } \left( \dsize
\prod_{i < \lambda^+} \lambda_i/J^{bd}_{\lambda^+} \right)$
\sn
\itemitem{ $(\varepsilon)$ }  for some regular $\kappa \le \lambda$, for any
$\mu' < \mu$ there is a tree $T$ with $\le \lambda$ nodes, $\kappa$ levels
and $|\text{Lim}_\kappa(T)| \ge \mu'$ (in fact e.g. $\kappa = \text{ Min}
\{ \kappa:2^\kappa \ge \mu\}$ is appropriate; without loss of generality
$T \subseteq {}^{\kappa >} \lambda$; we can get, of course, a tree $T$ with
cf$(\kappa)$ levels, $\le \lambda$ nodes and 
$|\text{Lim}_{\text{cf}(\kappa)}(T)| \ge \mu'$).
\endroster}
\sn
\item "{$(B)_\lambda$}"  for some $\mu,\chi^*$ we have: clauses $(\alpha) -
(\varepsilon)$ from above (so $2^\lambda < \chi^*$) and
{\roster
\itemitem{ $(\zeta)$ }  there is $\langle T_\zeta:\zeta < \chi^* \rangle$
such that: $T_\zeta \subseteq {}^{\lambda^+ >}2$ a tree, of cardinality
$\le \lambda^+$ and $2^{\lambda^+} = \dsize \bigcup_{\zeta < \chi^*}
\text{ Lim}_{\lambda^+}(T_\zeta)$ and $\chi^* < 2^{\lambda^+}$
\sn
\itemitem{ $(\eta)$ }  $2^\lambda \le \chi^* < \mu_{\text{wd}}
(\lambda^+,2^\lambda)$ (but $< \mu_{\text{wd}}(\lambda^+,2^\lambda)$ 
not used here, see \nl

$\quad$ \cite[Definition f.1(5)]{Sh:576})
\sn
\itemitem{ $(\theta)$ }  for some $\zeta < \chi^*$ we have Lim$_{\lambda^+}
(T_\zeta) \notin$ WDmTId$(\lambda^+)$
\sn
\itemitem{ $(\iota)$ }  if there is a normal $\lambda^{++}$-saturated ideal
on $\lambda^+$, e.g. the ideal WDmId$(\lambda^+)$ is, then 
$2^{\lambda^+} = \lambda^{++}$ (so as $2^\lambda < 2^{\lambda^+}$,
necessarily $2^\lambda = \lambda^+$)
\sn
\itemitem{ $(\kappa)$ }  cov$(\chi^*,\lambda^{++},\lambda^{++},\aleph_1) 
= \chi^*$, equivalently $\chi^* =$ \nl

$\quad \sup\{pp(\chi):\chi \le 2^\lambda,
\aleph_1 \le \text{cf}(\chi) \le \lambda^+ < \chi\}$ by \nl

$\quad$ \cite[Ch.II,5.4]{Sh:g}
\endroster}
\sn
\item "{$(C)_\lambda$}"  letting $\chi^* = 2^\lambda$ we have $(\zeta)$
(except ``$\chi^* < 2^{\lambda^+}$"), $(\eta),
(\theta),(\iota),(\kappa)$ of clause (B) and
{\roster
\itemitem{ $(\lambda)$ } for no $\mu \in (\lambda^+,2^\lambda]$ do we have 
cf$(\mu) \le \lambda^+$,pp$(\mu) >
2^\lambda$ equivalently cf$([2^\lambda]^{\lambda^+},\subseteq) = 2^\lambda$
hence $\mu_{wd}(\lambda^+,2^\lambda) = 2^{\lambda^+}$ except (maybe) when
$\lambda < \beth_\omega$ and there is ${\Cal A} \subseteq 
[\mu_{wd}(\lambda^+,2^\lambda)]^\lambda$ such that 
$A \ne B \in {\Cal A} \Rightarrow |A \cap B| = \aleph_0$.
\endroster}
\endroster
\enddemo
\bigskip

\remark{Remark}  Remember that

$$
\align
\text{cov}(\chi,\mu,\theta,\sigma) = \chi + \text{ Min} \bigl\{
|{\Cal P}|:&{\Cal P} \subseteq [\chi]^{< \mu} \text{ and every member of} \\
  &[\chi]^{< \theta} \text{ is included in the union of } < \sigma
\text{ members of } {\Cal P} \bigr\}.
\endalign
$$
\endremark
\bigskip

\demo{Proof}  This is related to \cite[II,5.11]{Sh:g}; we assume basic 
knowledge of pcf (or a readiness to believe).
Note that if $2^\lambda > \lambda^+$ then 
cf$([2^\lambda]^{\le \lambda^+},\subseteq) = 2^\lambda
\Leftrightarrow \text{ cov}(2^\lambda,\lambda^{++},\lambda^{++},\aleph_0) =
2^\lambda$ and cov$(2^\lambda,\lambda^{++},\lambda^{++},\aleph_0) \ge
\text{ cov}(2^\lambda,\lambda^{++},\lambda^{++},\theta) =
2^\lambda$ for $\theta \in [\aleph_0,\lambda]$.
\enddemo
\bigskip

\noindent
\underbar{Possibility 1}:  $\chi^* =: \text{ cov}(2^\lambda,\lambda^{++},
\lambda^{++},\aleph_1) = 2^\lambda$.

We shall show that case $(C)$ holds. \newline
Now by the definition of cov, clause $(\zeta)$ is obvious, as well as 
$(\kappa)$.  As on the one hand by \cite[AP,1.16 + 1.19]{Sh:f} we have
$\left( \mu_{wd}(\lambda^+,2^\lambda) \right)^{\aleph_0} = 2^{\lambda^+} > 
2^\lambda =
\chi^*$ and on the other hand $(\chi^*)^{\aleph_0} = (2^\lambda)^{\aleph_0} 
= 2^\lambda = \chi^*$ necessarily $\chi^* < \mu_{wd}(\lambda^+,2^\lambda)$ 
so clause $(\eta)$ follows; now clause $(\theta)$ follows from clause 
$(\zeta)$ as WDmTId$(\lambda^+)$ is $(2^\lambda)^+$-complete by 
\cite[1.2(5)]{Sh:576} and we have chosen $\chi^* = 2^\lambda$.  
Now if $2^{\lambda^+} > \lambda^{++}$, (so $2^{\lambda^+}
\ge \lambda^{+3}$), then for some $\zeta < \chi^*,T_\zeta$ is (a tree with
$\le \lambda^+$ nodes, $\lambda^+$ levels and) at least $\lambda^{+3} \,\, 
\lambda^+$-branches which is well known (see e.g. [J]) to imply ``no normal 
ideal on $\lambda^+$ is $\lambda^{++}$-saturated"; so we got clause $(\iota)$.
As for $(\lambda)$ the definition of $\chi^*$ and the assumption $\chi^* =
2^\lambda$ we have the first two phrases, as for 
$\mu_{wd}(\lambda^+,2^\lambda) =
2^{\lambda^+}$ by \cite[1.14 + 1.16]{Sh:f} there is ${\Cal A}$ as mentioned
in $(\lambda)$? and by \cite{Sh:460} we get $\lambda < \beth_\omega$.  The
``equivalently" holds as $(2^\lambda)^{\aleph_0} = 2^\lambda$.
\bigskip

\noindent
\underbar{Possibility 2}:  $\chi^* = \text{ cov}(2^\lambda,\lambda^{++},
\lambda^{++},\aleph_1) > 2^\lambda$.

Let $\mu = \text{ Min}\{\mu:\text{cf}(\mu) \le \lambda^+,\lambda^+ < \mu 
\le 2^\lambda$ and pp$(\mu) = \chi^*\}$, we know by \cite[II,5.4]{Sh:g} that 
$\mu$ exists and (by \cite[II,2.3(2)]{Sh:g}) clause $(\gamma)$ holds, also 
$2^\lambda < \text{ pp}(\mu) \le \mu^{\text{cf}(\mu)}$ hence 
cf$(\mu) = \lambda^+$.  So clauses $(\alpha),(\beta),(\gamma)$ hold 
(moreover, for clause $(\beta)$ use \cite[Ch.II,5.4(2)]{Sh:g}), 
and by $(\gamma)$ + \cite[VIII,\S1]{Sh:g} also clause $(\delta)$ holds.

For clause $(\varepsilon)$ let $\Upsilon = \text{ Min}\{\Upsilon:2^\Upsilon
\ge \mu\}$, clearly $\alpha < \Upsilon \Rightarrow 2^{|\alpha|} < \mu$ and
$\Upsilon \le \lambda$ (as $2^\lambda \ge \mu$) hence cf$(\Upsilon) \le
\Upsilon \le \lambda < \lambda^+ = \text{ cf}(\mu)$ hence 
$2^{<\Upsilon} < \mu$.  Now we shall first prove
\medskip
\roster
\item "{$(*)$}"  there \footnote{the less easy
point is when cf$(\Upsilon) = \aleph_0$,
otherwise we can get the conclusion differently (by \newline
\cite[II,5.4]{Sh:g}), so \scite{2.10}(A) suffice} 
is a tree with $\lambda^+$ nodes, cf$(\Upsilon)$ levels and 
$\ge \mu \,\Upsilon$-branches \newline
\smallskip
\noindent
why? otherwise we shall get contradiction to the claim \scite{2.10C}
below with $\sigma,\kappa,\theta_0,\theta_1,\mu,\chi$ there standing for 
cf$(\Upsilon),\lambda^+,\lambda^+,2^{< \Upsilon},\mu,(2^\lambda)^+$
here and $T^*$ defined below; let us check the conditions there:
\endroster
\medskip

\noindent
\underbar{Clause $(a)$}:  It says cf$(\Upsilon) < \lambda^+ = \text{ cf}
(\mu) \le \lambda^+ \le 2^{<\Upsilon} < \mu$ which is readily checked
except the inequality $\lambda^+ \le 2^{< \Upsilon}$ but if it fails we
immediately get more than required.
\medskip

\noindent
\underbar{Clause $(b)$}:  This is clause $(\gamma)$ of (A) which we have
proved.
\medskip

\noindent
\underbar{Clause $(c)$}:  The tree $T^*$ is 
$({}^{\Upsilon >}2,\triangleleft)$ restricted to an unbounded set of levels
of order type cf$(\Upsilon)$.
\mn
\underbar{Clause $(d)$}:  Let $\theta_2 =: \text{ cov}(2^{< \Upsilon},
\lambda^{++},\lambda^{++},\text{cf}(\Upsilon)^+)$.  \newline
So the statement we have to prove is pp$(\mu) \ge \chi = \text{ cf}(\chi)
> \theta^{\text{cf}(\Upsilon)}_2$.  Now pp$(\mu) \ge \chi$ holds by the
choice of $\mu$ and $\chi = \text{ cf}(\chi)$ as $\chi = (2^\lambda)^+$.
For the last inequality, by \cite[Ch.II,5.4]{Sh:g} and the choice of $\mu$,
as we have shown $2^{< \Upsilon} < \mu$ 
we know $\theta_2 < \mu$, but $\mu \le 2^\lambda$ so $\theta^{\text{cf}
(\Upsilon)}_2 \le (2^\lambda)^{\text{cf}(\Upsilon)} \le (2^\lambda)^\Upsilon 
\le (2^\lambda)^\lambda = 2^\lambda < \chi$.
\medskip

So we have verified clauses $(a)-(d)$ of \scite{2.10C} hence its conclusion
holds, but this gives $(*)$, i.e. the desired conclusion in clause 
$(\varepsilon)$ of Case A in \scite{2.10}; well not exactly, it gives 
only $|T^*| \le \lambda^+$, so $T^* = \dsize
\bigcup_{i < \lambda^+} T_i,T_i$ increases continuously with each $T_i$
of cardinality $\le \lambda$, so for every $\mu' < \mu$ for some $i$ 
we have 
$|\text{Lim}_{\text{cf}(\Upsilon)}(T_i)| \ge \mu'$,
so we have clauses $(\alpha)-(\varepsilon)$ there. 
\bn
\underbar{Subpossibility 2a}:  $\chi^* < 2^{\lambda^+}$.  \nl

We shall prove $(B)_\lambda$, so we are left with 
proving $(\zeta)-(\kappa)$ when $\chi^* < 2^{\lambda^+}$.
By the choice of $\chi^*$, easily clause $(\zeta)$ (in Case B of
\scite{2.10}) holds.  In clause $(\eta),``2^\lambda < \chi^*"$ holds 
as we are in possibility 2.
\medskip

Also as pp$(\mu) = \chi^*$ by the choice of $\mu$ necessarily 
(by transitivity of pcf, i.e. \cite[Ch.II,2.3(2)]{Sh:g}) cf$(\chi^*) > 
\lambda^+$ but $\mu > \lambda^+$.  
Easily $\chi \le \chi^* \wedge$ cf$(\chi) \le \lambda^+
\Rightarrow$ pp$(\chi) \le \chi^*$ hence cov$(\chi^*,\lambda^{++},
\lambda^{++},\aleph_1) = \chi^*$ by \cite[Ch.II,5.4]{Sh:g}, which gives clause
$(\kappa)$, so as $(\lambda^+)^{\aleph_0} \le 2^\lambda < \chi^*$ 
certainly there is no family of $>
\chi^*$ subsets of $\chi^*$ each of cardinality $\lambda^+$ with pairwise
finite intersections, hence (by \cite[Ch.XIV,\S1]{Sh:b} or see \relax
\cite[1.2(1)]{Sh:576} or \cite[AP,1.16]{Sh:f}) we have 
$\chi^* < \mu_{wd}(\lambda^+,2^\lambda)$ thus 
completing the proof of $(\eta)$.

Now clause $(\theta)$ follows by $(\zeta) + (\eta)$ by \cite[1.2(5)]{Sh:576}.
Also if
$2^{\lambda^+} \ne \lambda^{++}$ then $2^{\lambda^+} \ge \lambda^{+3}$ so by
clause $(\zeta)$ (as $\chi^* < 2^{\lambda^+}$), for some $\zeta,
|\text{Lim}_{\lambda^+}(T_\zeta)| \ge \lambda^{+3}$, which is well known
to imply no normal ideal on $\lambda^+$ is $\lambda^{++}$-saturated; i.e.
clause $(\iota)$.  So we have proved that case $(B)_\lambda$ holds.
\bigskip

\noindent
\underbar{Subpossibility 2b}:  $\chi^* = 2^{\lambda^+}$.

We have proved that case $(A)_\lambda$ holds, as we already defined $\mu$ 
and $\chi^*$ and proved
$(\alpha),(\beta),(\gamma),(\delta),(\varepsilon)$ we are done.
\bigskip

\noindent
Still we depend on \scite{2.10C} below which in turn depends on \scite{2.10B}.
\proclaim{\stag{2.10B} Claim}  Assume
\medskip
\roster
\item "{$(a)$}"  $\sigma < \kappa = \text{ cf}(\mu) \le \theta_0 \le \theta_1
< \mu \le \theta^\sigma_1$
\sn
\item "{$(b)$}"  $(\forall \mu')[\theta_0 < \mu' < \mu \and \text{ cf}
(\mu') \le \kappa \Rightarrow \text{ pp}(\mu') < \mu]$
\sn
\item "{$(c)$}"  $\theta_2 = \theta_1 +\text{ cov}
(\theta_1,\theta^+_0,\kappa^+,\sigma^+)$ (by clause (b) and
\cite[Ch.II,5.4]{Sh:g} we know that it is $< \mu$)
\sn
\item "{$(d)$}"  pp$(\mu) \ge \chi = \text{ cf}(\chi) > \theta^\sigma_2
\,\, (\ge \theta^\sigma_1 \ge \mu$).
\endroster
\medskip

\underbar{Then} $\theta^\sigma_0 \ge \mu$.
\endproclaim
\bigskip

\remark{Remark}  In fact $\theta^\sigma_2 \ge \text{ cov}(\theta_1,
\theta^+_0,\kappa^+,2)$.
\endremark
\bigskip

\demo{Proof}  Assume toward contradiction $\theta^\sigma_0 < \mu$.  By
\cite[Ch.II,2.3(2)]{Sh:g} and clause (b) of the assumption we have 
sup$\{\text{pp}(\mu'):\theta^+_0 \le \mu' \le \theta^\sigma_0$ and 
cf$(\mu') \le \kappa\} < \mu$ hence by \cite[Ch.II,5.4]{Sh:g} it follows that 
$\kappa^* =: \text{ cov}(\theta^\sigma_0,\theta^+_0,\kappa^+,\sigma^+) < \mu$.
We can by assumption (b) + (d) and \cite[Ch.II,3.5]{Sh:g} + 
\cite[Ch.VIII,\S1]{Sh:g} find $T \subseteq {}^{\kappa \ge}\mu$ a tree with
$\le \mu$ nodes, $|\text{Lim}_\kappa(T)| \ge \chi$, 
(if $\chi = \text{ pp}(\mu)$, the supremum in the definition of pp$(\mu)$ is 
obtained by \cite[II,5.4(2)]{Sh:g}).  Moreover, by the construction there is
$\Xi \subseteq \text{ Lim}_\kappa(T),|\Xi| = \chi$ such that $\Xi' \subseteq
\Xi \and |\Xi'| \ge \chi \Rightarrow |\{ \eta \restriction \alpha:\alpha <
\kappa,\eta \in \Xi'\}| = \mu$.  By renaming (and also by the construction), 
\wilog 
\mr
\item "{$\otimes$}"  if $\eta_0 \char 94 \langle \alpha_0 \rangle \ne
\eta_1 \char 94 \langle \alpha_1 \rangle$ belongs to $T$ then
$\alpha_0 \ne \alpha_1$.
\ermn
So let $\eta_i \in \text{ Lim}_\kappa(T)$ for $i < \chi$ be pairwise
distinct, listing $\Xi$. \newline
As $\mu \le \theta^\sigma_1$ there is a sequence
$\bar F = \langle F_\varepsilon:
\varepsilon < \sigma \rangle$ satisfying: $F_\varepsilon$ a function from
$\mu$ to $\theta_1$ such that $\alpha < \beta < \mu \Rightarrow (\exists
\varepsilon < \theta)F_\varepsilon(\alpha) \ne F_\varepsilon(\beta)$.

Let $w_{i,\varepsilon} = \{F_\varepsilon(\eta_i(\alpha)):\alpha < \kappa\}$,
so $w_{i,\varepsilon} \in [\theta_1]^\kappa$.  By assumption (c) we
have \newline
$\theta_2 = \theta_1 +
\text{ cov}(\theta_1,\theta^+_0,\kappa^+,\sigma^+)$ so there is ${\Cal P}
\subseteq [\theta_1]^{\theta_0},\theta_2 = |{\Cal P}|$ such that: any 
$w \in [\theta_1]^\kappa$ is included in a union of $\le \sigma$ members of
${\Cal P}$.  So we can find $X_{i,\varepsilon,\zeta} \in {\Cal P}$ for
$\zeta < \sigma$ such that $w_{i,\varepsilon} \subseteq \dsize \bigcup
_{\zeta < \sigma} X_{i,\varepsilon,\zeta}$.  So $\dsize \bigcup_{\varepsilon
< \sigma} w_{i,\varepsilon} \subseteq Y_i =: \dsize \bigcup_{\zeta,
\varepsilon < \sigma} X_{i,\varepsilon,\zeta}$.  Let
${\Cal P}^* = \{ \dsize \bigcup_{\varepsilon < \sigma} X_\varepsilon:
X_\varepsilon \in {\Cal P}$ for $\varepsilon < \sigma\}$, so ${\Cal P}^*$ is
a family of $\le |{\Cal P}|^\sigma \le \theta^\sigma_2$ sets and $i < \chi
\Rightarrow Y_i \in {\Cal P}^*$.
\medskip

For each $Y \in {\Cal P}^*$ let

$$
Z_Y = \{ \alpha < \mu:(\forall \varepsilon < \sigma)(F_\varepsilon(\alpha)
\in Y)\}
$$
\medskip

\noindent
clearly $Y = Y_i \Rightarrow \text{ Rang}(\eta_i) \subseteq Z_Y$, also
$|Y| \le \theta_0$ hence $|Z_Y| \le \theta^\sigma_0 < \mu$ hence there is a
family $Q_Y$ of cardinality $\kappa^* =: \text{ cov}(\theta^\sigma_0,
\theta^+_0,\kappa^+,\sigma^+) < \mu$ whose members are subsets of $Z_Y$ each 
of cardinality $\le \theta_0$
such that any $X \in [Z_Y]^{\le \kappa}$ is included in the union of
$\le \sigma$ of them.  For each $Y \in {\Cal P}^*$ and $W \in Q_Y$ let
$T'_W = \{ \eta \in T:\text{for some } \alpha < \kappa \text{ we have: }
\alpha + 1 = \ell g(\eta)$ and $\eta(\alpha) \in W\}$ and 
$T_W = \{\eta \in T'_W:(\exists \nu)(\eta \triangleleft \nu \in T_W)\}$.

So by $\otimes$ above we have: $T'_W$, hence $T_W$ is a set of 
$\le |W| + \kappa \le \theta_0$ nodes in $T$, $\triangleleft$-downward 
closed.  Also
\medskip
\roster
\item "{$(*)$}"  $| \dsize \bigcup_{Y \in {\Cal P}^*} Q_Y| \le
|{\Cal P}^*| \times \underset{Y \in {\Cal P}^*} {}\to \sup |Z_Y|$ \newline
$\le \theta^\sigma_2 \times \kappa^\sigma \le \theta^\sigma_2 +
\text{ cov}(\theta^\sigma_0,\theta^+_0,\theta^+_0,\sigma^+) < \mu$.
\endroster
\medskip

However, for every 
$i < \chi,Y_i \in {\Cal P}^*$ and Rang$(\eta_i) \in [Y_i]^{\le \kappa}$ so
for some \newline
$W \in Q_{Y_i},(\exists^\kappa \alpha < \kappa)[\eta_i(\alpha) 
\subseteq W_i]$ hence $\eta _i \in \text{ Lim}_\kappa(T_W)$.

By assumption (a) and $(*)$ above for some $W \in \dsize \bigcup
_{Y \in {\Cal P}^*} Z_Y$ we have

$$
|\{i < \chi:\eta_i \in \text{ Lim}_\kappa(T_W)\}| = \chi.
$$
\medskip

\noindent
$T_W$ is (essentially) a tree and contradict the choice of $\Xi = \{\eta_i:
i < \chi\}$. \nl
(We could have instead using $\kappa^*,Q_Y$ to fix $Y_i = Y$ as 
$|{\Cal P}^*| < \chi = \text{ cf}(\chi)$.) \hfill$\square_{\scite{2.10B}}$
\enddemo
\bigskip

\proclaim{\stag{2.10C} Claim}  Assume
\roster
\item "{$(a)$}"  $\sigma < \kappa = \text{ cf}(\mu) \le \theta_0 \le
\theta_1 < \mu$
\sn
\item "{$(b)$}"  $(\forall \mu')[\theta_0 < \mu' < \mu \and \text{ cf}
(\mu') \le \kappa \rightarrow \text{ pp}(\mu') < \mu]$
\sn
\item "{$(c)$}"  $T^*$ is a tree with $\le \theta_1$ nodes, $\sigma$ levels 
and $\ge \mu \,\,\sigma$-branches
\sn
\item "{$(d)$}"  pp$(\mu) \ge \chi = \text{ cf}(\chi) > \theta^\sigma_2$
where $\theta_2 = \text{ cov}(\theta_1,\theta^+_0,\kappa^+,\sigma^+)$.
\endroster
\medskip

\noindent
\underbar{Then} for some subtree $Y \subseteq T^*,|Y| \le \theta_0$ and
$|\text{Lim}_\sigma(Y)| \ge \mu$ (enough $\ge \mu'$ for any given $\mu' <
\mu$).
\endproclaim
\bigskip

\demo{Proof}  Let $T,\Xi = \{ \eta_i:i < \chi\}$ be as in the proof of the
previous claim.  Let \nl
$\{\nu_\zeta:\zeta < \mu\}$ list $\mu$ distinct
$\sigma$-branches of $T^*$ (see clause (c)).  
Without loss of generality the set of
nodes of $T^*$ is $\theta_1$.  Choose for each $\varepsilon < \sigma$ the 
function $F_\varepsilon:\mu \rightarrow \theta_1$ by $F_\varepsilon(\gamma) =
\nu_\gamma(\varepsilon)$.  Define $w_{i,\varepsilon},{\Cal P},
X_{i,\varepsilon,\zeta},Y_i,{\Cal P}^*,Z_Y$ as in the proof of \scite{2.10B}.
But for $Y \in {\Cal P}^*$ we change the choice of $Z_Y$, first

$$
Y' = \{ \beta < \theta_1:\text{for some } \alpha \in Y, \text{ we have }
\beta <_{T^*} \alpha\}
$$
\medskip

\noindent
So $|Y'| \le \sigma + |Y|$ and let

$$
Z_Y = \{\alpha < \mu:(\forall \varepsilon < \sigma)(F_\varepsilon(\alpha)
\in Y')\}.
$$
\medskip

\noindent
We continue as in the proof of \scite{2.10B}. \hfill$\square_{\scite{2.10C}},
\,\,\square_{\scite{2.10}}$
\enddemo
\bigskip

\remark{\stag{2.10D} Remark}  1) We could have used in \scite{2.10B}, 
\scite{2.10C}, \nl
$\theta_2 = \text{ cov}(\theta^\sigma_0,\theta^+_0,\kappa^+,
\kappa^+)$ instead cov$(\theta^\sigma_0,\theta^+_0,\kappa^+,\sigma^+)$ and
similarly in the proof of \scite{2.10}. \nl
2) We can also play with assumption (b) as \scite{2.10B}, \scite{2.10C}.
\endremark
\bigskip

\noindent
It may be useful to note
\demo{\stag{2.11} Fact}  If $T \subseteq {}^{\lambda^+ >}2$ is a tree, 
$|T| \le \lambda^+$ and $\lambda \ge \beth_\omega$ \underbar{then} for every
regular $\kappa < \beth_\omega$ large enough, we can find $\langle Y_\delta:
\delta < \lambda^+,\text{cf}(\delta) = \kappa \rangle,|Y_\delta| \le
\lambda$ such that: \newline
for every $\eta \in \text{ Lim}_{\lambda^+}(T)$ for a club
of $\delta < \lambda^+$ we have \newline
cf$(\delta) = \kappa \Rightarrow \eta \restriction \delta \in Y_\delta$.
\enddemo
\bigskip

\demo{\stag{2.12} Fact}  Assume ${\frak K}$ is an abstract elementary class 
with amalgamation in $\lambda$, and above $(M,N,a) \in K^3_\lambda$ there 
is no minimal pair. \newline
1)   Assume $T$ is a tree with $\delta < \lambda^+$ levels and $\le \lambda$ 
nodes.  \ub{Then} we can find $(M^*,N_\eta,a) \in K^3_\lambda$ above
$(M,N,a)$ for $\eta \in \text{ Lim}_\delta(T)$ such that tp$(a,M^*,N_\eta)$
for $\eta \in \text{ Lim}_\delta(T)$ are pairwise distinct.  We can add 
``$(M^*,N_\eta,a)$ is reduced". \newline
2) If $M \in K_\lambda$ is universal then ${\Cal S}(M) \ge
\sup\{\text{Lim}_\delta(T):T$ a tree with $\le \lambda$ nodes and $\delta$
levels$\}$.
\enddemo
\bigskip

\demo{Proof}  1) Straight (or see the proof of \scite{2.16}(1)). \nl
2) As for any $N \in K_\lambda$ there is a model $N' \le_{\frak K} M$
isomorphic to $N$, now $p \mapsto p \restriction N'$ is a function from
${\Cal S}(M)$ onto ${\Cal S}(N')$ by \cite[\S2]{Sh:576} 
hence $|{\Cal S}(M)| \ge |{\Cal S}(N')| =
|{\Cal S}(N)|$.  Now use part (1).  \hfill$\square_{\scite{2.12}}$
\enddemo
\bigskip

\noindent
Recall from \cite[3.13(3)]{Sh:576}.
\proclaim{\stag{2.14} Claim}  1) Assume
\medskip
\roster
\item "{$(a)$}"  $2^\lambda < 2^{\lambda^+}$ and Case A or B of 
Fact \scite{2.10} holds for $\mu,\chi^*$ (or just the conclusion there)
\sn
\item "{$(b)$}"  ${\frak K}$ is an abstract elementary class with
$LS({\frak K}) \le \lambda$
\sn
\item "{$(c)$}"  $K_{\lambda^+} \ne 0$
\sn
\item "{$(d)$}"  ${\frak K}$ has amalgamation in $\lambda$
\sn
\item "{$(e)$}"  in $K^3_\lambda$, the minimal triples are not dense.
\endroster
\medskip

\noindent
\underbar{Then}
\medskip
\roster  
\item "{$(*)_1$}"  for any regular $\chi \le \mu$ we have:
\sn
\item "{$(*)^1_\chi$}"  there is $M \in K_\lambda,|{\Cal S}(M)| \ge \chi$.
\endroster
\medskip

\noindent
2) If in part (1) we strengthen clause $(d)$ to $(d)^+$, \ub{then} 
we get $(*)^+_1$ where:
\mr
\item "{$(d)^+$}"  ${\frak K}$ has amalgamation in $\lambda$ and a 
universal member in $\lambda$
\sn
{\roster
\itemitem{ $(*)^+_1$ }  for some $M \in K_\lambda$ we have 
$|{\Cal S}(M)| \ge \mu$.
\endroster}
\ermn
3)  Assume (a), (b), (c), (e) of part (1) and (d)$^+$ of part (2) 
\underbar{then}:
\medskip
\roster
\item "{$(*)_2$}"  $I(\lambda^+,K) \ge \chi^*$ and 
if $(2^\lambda)^+ < \chi^*$ then $IE(\lambda^+,{\frak K}) \ge \chi^*$.
\endroster
\medskip

\noindent
4) If in clause $(a)$ of part (1) we restrict ourselves to Case A of
\scite{2.10}, \underbar{then} \nl
$\chi^* = 2^{\lambda^+}$ so in part (3) we get
\medskip
\roster
\item "{$(*)^+_2$}"  $I(\lambda^+,K) = 2^{\lambda^+}$ and
$(2^\lambda)^+ < 2^{\lambda^+} \Rightarrow IE(\lambda^+,K) \ge 2^{\lambda^+}$.
\endroster 
\endproclaim
\bigskip

\remark{\stag{2.15} Remark}  1) We can restrict clause (b) to 
$K_\lambda$, interpreting in (c) + (e), $K_{\lambda^+}$ as 
$\{ \dsize \bigcup_{i < \lambda^+} M_i:M_i \in K_\lambda\}$ is 
$<_{\frak K}$-increasing (strictly and) continuous, but see \S0, mainly
\scite{0.31}.
\newline
2) Part (3) of \scite{2.14} is like \cite[Ch.II,4.10E]{Sh:g}, Kojman Shelah 
\cite[\S2]{KjSh:409}. \newline
3) We can apply this to $\lambda^+$ standing for $\lambda$ . \newline
4) We can state the part of (A) of \scite{2.10} used (and can replace
$2^{\lambda^+}$ by smaller cardinals). \newline
5) We can replace $\lambda^+$ by a weakly inaccessible cardinal with
suitable changes.
\endremark
\bigskip

\demo{Proof}   1) Note that $\mu$ is singular (as by clause $(\alpha)$ of
(A) of \scite{2.10}, cf$(\mu) = \lambda^+ < \mu$).
By \scite{2.12}(1) it suffices for each $\mu' < \mu$ to have 
$\delta < \lambda^+$ and a tree with $\le \lambda$ nodes and $\ge \mu' \,\,
\delta$-branches.  They exist by clause $(\varepsilon)$ of (A) of 
\scite{2.10}.  \newline
2) Similarly using \scite{2.12}(2). \newline
3) So assume $M^* \in K_\lambda$ and ${\Cal S}(M^*)$ has cardinality 
$\ge \mu$, let $p_\eta \in {\Cal S}(M^*)$ for $\eta \in Z$ be pairwise 
distinct, $|Z| \ge \mu$ and let $M^* \le_{\frak K} N_\eta \in K_\lambda,
p_\eta = \text{ tp}(a_\eta,M^*,N_\eta)$. \newline
Now for every $X \in [Z]^{\lambda^+}$, as ${\frak K}$ has
amalgamation in $\lambda$ there is $M_X \in K_{\lambda^+}$ such that
$M^* \le_{\frak K} M_X$ and $\eta \in X \Rightarrow N^*_\eta$ is embeddable 
into $M_X$ over $M^*$ (hence $p_\eta$ is realized in $M_X$).  
Let $Y[X] = \{\eta \in Z:p_\eta \text{ is realized in } M_X\}$.  So 
$X \subseteq Y[X] \in [Z]^{\lambda^+}$, so $\{Y[X]:X \in [Z]^{\lambda^+}\}$ 
is a cofinal subset of $[Z]^{\lambda^+}$, hence 

$$
\align
|\{(M_X,c)_{c \in M^*}/\cong :&\,X \in [Z]^{\lambda^+}\}| \ge \\
  &\,|\{Y[X]:X \in [Z]^{\lambda^+}\}| \ge \text{ cf}([Z]^{\lambda^+},
\subseteq) \ge \\
  &\,\text{cf}([\mu]^{\lambda^+},\subseteq) \ge pp(\mu) = \chi^*.
\endalign
$$
\medskip

As $2^\lambda < \chi^*$ also $|\{M_X/\cong :X \in [Z]^{\lambda^+}\}| 
\ge \chi^*$ (clear or see \cite[Ch.VIII,1.2]{Sh:a} because $\|M_X\| =
\lambda^+,\|M^*\| = \lambda$ and $(\lambda^+)^\lambda < \mu$) but 
$I(\lambda^+,K)$ is $\ge$ than the former.

Lastly we shall prove $(2^\lambda)^+ < 2^{\lambda^+} \Rightarrow IE(\lambda^+,
K) \ge \chi^*$ (so the reader may skip this, sufficing himself with the 
estimate on $I(\lambda^+,K)$). \newline
For each $X \in [\mu]^{\lambda^+}$, let $\bold F_X = \{f:f \text{ a }
\le_{\frak K} \text{-embedding of } M^* \text{ into } M_X\}$, and for $f \in
F_X$ let

$$
{\Cal Z}_{X,f} = \bigl\{X_1 \in [Z]^{\lambda^+}:\text{there is a }
\le_{\frak K}\text{-embedding of } M_{X_1} \text{ into } M_X 
\text{ extending } f \bigr\},
$$
\medskip

\noindent
and let ${\Cal S}_{X,f} = \{p \in {\Cal S}(M^*):f(p) \text{ is realized in }
M_X\}$, so ${\Cal Z}_{X,f} \subseteq \{X_1:X_1 \subseteq {\Cal S}_{X,f}\}$ and
$|{\Cal S}_{X,f}| \le \lambda^+$.

Now the result follows from the the fact \scite{2.15A} below. \nl
4) Should be clear from the proof of part (3).
\enddemo
\bigskip

\noindent
\underbar{\stag{2.15A} Fact}:  Assume:
\medskip
\roster
\item "{$(a)$}"  cf$(\mu) \le \kappa < \mu$, pp$_\kappa(\mu) = \chi^*$, 
moreover pp$^+_\kappa(\mu) = (\chi^*)^+$ and $\kappa^+ < \theta < \chi^*$
\sn
\item "{$(b)$}"  $\bold F$ is a function, with domain $[\mu]^\kappa$, such
that: for $a \in [\mu]^\kappa,\bold F(a)$ is a family of $< \theta$ members of
$[\mu]^\kappa$
\sn
\item "{$(c)$}"  $F$ is a function with domain $[\mu]^\kappa$ such that
$$
a \in [\mu]^\kappa \Rightarrow a \subseteq F(a) \in \bold F(a).
$$   
\endroster
\medskip

\noindent
\underbar{Then} we can find pairwise distinct $a_i \in [\mu]^\kappa$ for
$i < \chi^*$ such that ${\Cal I} = \{a_i:i < \chi\}$ is 
$(F,\bold F)$-independent which means
\medskip
\roster
\item "{$(*)_{F,\bold F,{\Cal I}}$}"  $\qquad$ if 
$a \ne b \and a \in {\Cal I} \and
b \in {\Cal I} \and c \in \bold F(a) \Rightarrow \neg(F(b) \subseteq c)$.
\endroster
\bigskip

\remark{Remark}  1) Similar to Hajnal's free subset theorem \cite{Ha61}.
Without loss of generality $F(a) = F(b)$.  No loss in assuming $F(a) = a$.
\nl
2) Note that we can let $F(a) = a$. \nl
3) Note that if $\lambda = \text{ cf}([\mu]^\kappa,\subseteq)$ then for some
$\bold F$ as in the Fact
\mr
\item "{$(*)$}"  if $a_i \in [\mu]^\kappa$ for $i < \lambda^+$ are pairwise
distinct then no pair $\{a_i,a_j\}$ is $(\bold F,F)$-independent \nl
[why?  let ${\Cal P} \subseteq [\mu]^\kappa$ be cofinal (under $\subseteq$)
of cardinality $\lambda$, and let $\bold F$ be such that \nl
$\bold F(a) \subseteq \{b \in [\mu]^\kappa:a \subseteq \kappa \text{ and }
b \in {\Cal P}\}$; \nl
clearly there is such $\bold F$.  Now clearly
{\roster
\itemitem{$ (*)_1$ }  if 
$a \ne b$ are from $[\mu]^\kappa$ and $\bold F(a) \cap
\bold F(b) \ne \emptyset$ then $\{a,b\}$ is not $(\bold F,F)$-independent.
\nl
Also if $\mu_1 \le \mu$, cf$(\mu_1) \le \kappa \le \kappa + \theta < \mu_1$
and pp$_\kappa(\mu_1) \le$? then by \cite[Ch.II,2.3]{Sh:g} the Fact for $\mu_1$
implies the one for $\mu$.]
\endroster}
\endroster
\endremark  
\bigskip

\demo{Proof}  We can define $g:[\mu]^{\le \kappa} \rightarrow [\mu]^\kappa,
F',\bold F'$ as follows:

$$
g(a) = \{\kappa + \alpha:\alpha \in a\} \cup \{\alpha:\alpha < \kappa\}
$$

$$
\bold F'(a) = \{\{\alpha:\kappa + \alpha \in b\}:b \in \bold F'(g(a))\}
$$

$$
F(a) = \{\alpha:\kappa + \alpha \in F(g(a))\}.
$$
\mn
Now $\bold F,F$ are as above only replacing everywhere $[\mu]^\kappa$ by
$[\kappa]^\kappa$, and if ${\Cal I} = \{a_i:i < \chi\} \subseteq
[\mu]^{\le \kappa}$ with no repetitions satisfying $(*)_{F',\bold F',I}$
then ${\Cal I}' = \{g(a_i):i < \chi\}$ is with no repetitions and
$(*)_{F,\bold F,{\Cal I}'}$.
\sn
So we conclude that we can replace $[\mu]^{\le \kappa}$ by 
$[\mu]^{|{\frak a}|}$ where ${\frak a}$ is chosen below, i.e. we can replace
$\kappa$ by $|{\frak a}|$.
\sn
Without loss of generality $\kappa^{++} < \theta$.  

Assume $\theta < \chi = \text{ cf}(\Pi {\frak a}/J)$ where ${\frak a}
\subseteq \mu \cap \text{ Reg} \backslash \kappa^+,|{\frak a}| \le \kappa,
\sup({\frak a}) = \mu,J^{bd}_{\frak a} \subseteq J$ and for simplicity
$\chi = \text{ max pcf}({\frak a})$ so as explained above \wilog \nl
$|{\frak a}| < \kappa$ and let $f = \langle f_\alpha:\alpha <
\chi \rangle$ be a sequence of members of $\Pi{\frak a},<_J$-increasing, and
cofinal in $(\Pi{\frak a},<_J)$, so, of course, $\chi \le \chi^*$.  
Without loss of generality $f_\alpha(\tau)
> \sup({\frak a} \cap \tau)$ for $\tau \in {\frak a}$.  Also for every
$a \in [\mu]^\kappa$, define $f_a \in \Pi{\frak a}$ by $f_a(\tau) = \sup(a
\cap \tau)$ for $\tau \in {\frak a}$ so for some $\zeta(a) < \chi$ we have
$f_a <_J f_{\zeta(a)}$ (as $\langle f_\alpha:\alpha < \chi \rangle$ is
cofinal in $(\Pi{\frak a},<_J)$).  So for each $a \in [\mu]^\kappa$, as
$|\bold F(a)| < \theta < \chi = \text{ cf}(\chi)$ clearly 
$\xi(a) = \sup\{\zeta(b):b \in
\bold F(a)\}$ is $< \chi$, and clearly $(\forall b \in \bold F(a))
[f_b <_J f_{\xi(a)}]$.  So \newline
$C = \{\gamma < \chi:\text{for every } \beta < \gamma,\xi(\kappa \cup
\text{Rang }f_\beta) < \gamma\}$ is a club of $\chi$.

For each $\alpha < \chi$, Rang$(f_\alpha) \cup \kappa \in [\mu]^\kappa$, hence
$\bold F(\text{Rang}(f_\alpha) \cup \kappa)$ has cardinality $< \theta$, but 
$\theta <  \chi = \text{ cf}(\chi)$ hence for some $\theta_1 < \theta$; we 
have $\theta_1 > \kappa^{++}$ and $\chi = \sup\{\alpha < \chi:|\bold F
(\text{Rang}(f_\alpha) \cup \kappa)| \le \theta_1\}$, so without loss 
of generality $\alpha < \chi \Rightarrow \theta_1 \ge |\bold F(\text{Rang}
(f_\alpha) \cup \kappa)|$.
\medskip

As $\kappa^+ < \theta_1$, there is a stationary $S \subseteq \{\delta <
\theta^+_1:\text{cf}(\delta) = \kappa^+\}$ which is in $I[\theta^+_1]$, 
by \cite[\S1]{Sh:420} and
let $\langle d_i:i < \theta^+_1 \rangle$ witness it, so otp$(d_i) \le
\kappa^+,d_i \subseteq i,[j \in d_i \Rightarrow d_j = d_i \cap i]$ and
$i \in S \Rightarrow i = \sup(d_i)$, and for simplicity: for every club $E$ of
$\theta^+_1$ for stationarily many $\delta \in S$ we have $(\forall \alpha
\in d_\delta)[(\exists \beta \in E)(\sup(\alpha \cap d_\delta) < \beta <
\alpha)]$, exists by \cite[\S1]{Sh:420}.  
Now try to choose by induction on $i < \theta^+_1$, a triple
$(g_i,\alpha_i,w_i)$ such that:
\medskip
\roster
\item "{$(a)$}"   $g_i \in \Pi{\frak a}$
\smallskip
\noindent
\item "{$(b)$}"  $j < i \Rightarrow g_j <_J g_i$
\smallskip
\noindent
\item "{$(c)$}"  $(\forall \tau \in {\frak a})
(\underset{j \in d_i} {}\to \sup g_j(\tau) < g_i(\tau))$
\smallskip
\noindent
\item "{$(d)$}"  $\alpha_i < \chi,\alpha_i > \sup(\dsize \bigcup_{j<i} w_j)$
\smallskip
\noindent
\item "{$(e)$}"  $j < i \Rightarrow \alpha_j < \alpha_i$
\smallskip
\noindent
\item "{$(f)$}"  $g_i <_J f_{\alpha_i}$
\smallskip
\noindent
\item "{$(g)$}"  $\beta \in \dsize \bigcup_{j<i} w_j \Rightarrow \xi(\beta)
< \alpha_i \and f_\beta <_J g_i$
\smallskip
\noindent
\item "{$(h)$}"  $w_i$ is a maximal subset of $(\alpha_i,\chi)$ satisfying
{\roster
\itemitem{ $(*)$ }  $\beta \in w_i \and \gamma \in w_i \and \beta \ne \gamma
\and a \in \bold F(\kappa \cup \text{ Rang}(f_\beta)) 
\Rightarrow \neg(F(\text{Rang } f_\gamma) \subseteq a)$
\endroster}
or just
{\roster
\itemitem{ $(*)^+$ }   $\beta \in w_i \and \gamma \in w_i \and \beta \ne
\gamma \and a \in \bold F(\kappa \cup \text{ Rang}(f_\beta)) \Rightarrow$ \nl
$\{\tau \in {\frak a}:f_\gamma (\tau) \in a \} \in J$.
\endroster}
[note that really 
\item "{$\otimes$}"  if $w \subseteq (\alpha_i,\chi)$ satisfies $(*)^+$ then
it satisfies $(*)$
\sn
why?  let us check $(*)$, so let $\beta \in w,\gamma \in w,\beta \ne \gamma$
and $a \in \bold F(\text{Rang}(f_\beta))$; by $(*)^+$ we know that
${\frak a}' = \{\tau \in {\frak a}:f_\gamma(\tau) \in a\} \in J,J$ is a
proper ideal on ${\frak a}$ clearly for some $\tau \in {\frak a}$ we have
$\tau \notin {\frak a}'$, hence $f_\gamma(\tau) \notin a$ but $f_\gamma(\tau)
\in \text{ Rang}(f_\gamma) \subseteq F(\text{Rang}(f_\gamma))$ hence
$f_\gamma(\tau) \in F(\text{Rang}(f_\gamma)) \backslash a$ so $\neg(F
(\text{Rang}(f_\gamma) \subseteq a)$, as required.]
\endroster
\medskip

\noindent
We claim that we cannot carry the induction because if we succeed, then as
cf$(\chi) = \chi > \theta \ge \theta^+_1$ there is $\alpha$ such that
$\dsize \bigcup_{i < \theta^+_1} \alpha_i < \alpha < \chi$ and let
$\bold F(\kappa \cup \text{ Rang}(f_\alpha)) = 
\{a^\alpha_\zeta:\zeta < \theta_1\}$ (possible as
$1 \le |\bold F(\kappa \cup \text{ Rang}(f_\alpha))| \le \theta_1)$.  
Now for each $i < \theta^+_1$, by the choice of $w_i$ clearly $w_i \cup 
\{\alpha\}$ does not
satisfy the demand in clause $(h)$, so as $\beta \in w_i \Rightarrow
\xi(\beta) < \alpha_{i+1} < \alpha$, necessarily for some $\beta_i \in
w_i$ and $\zeta_i < \theta_1$ we have

$$
{\frak a}_i = \{\tau \in {\frak a}:f_{\beta_i}(\tau) \in a^\alpha_{\zeta_i}\}
\notin J.
$$
\mn
[why use the ideal?  In order to show that ${\frak b}_\varepsilon \ne
\emptyset$.]
Now cf$(\theta^+_1) = \theta^+_1 > \theta_1$, for some $\zeta(*) <
\theta^+_1$ we have $A = \{i:\zeta_i = \zeta(*)\}$ is unbounded in 
$\theta^+_1$.  Hence $E = \{\alpha < \theta^+_1:\alpha$ a limit ordinal and
$A \cap \alpha$ is unbounded in $\alpha\}$ is a club of $\theta^+_1$.  So for
some $\delta \in S$ we have $\delta = \sup(A \cap \delta)$, moreover if
$d_\delta = \{\alpha_\varepsilon:\varepsilon < \kappa^+\}$ (increasing) then
$(\forall \varepsilon)[E \cap (\underset{\zeta < \varepsilon} {}\to \sup
\,\alpha_\zeta,\alpha_\varepsilon) \ne \emptyset]$ hence we can find
$i(\delta,\varepsilon) \in (\underset{\zeta < \varepsilon} {}\to \sup \,
\alpha_\zeta,\alpha_\varepsilon) \cap A$ for each $\varepsilon < \kappa^+$.
\medskip

\noindent
Clearly for each $\varepsilon < \kappa^+$

$$
\align
{\frak b}_\varepsilon = \bigl\{ \tau \in {\frak a}:g_{i(\delta,\varepsilon)}
(\tau) &< f_{\alpha_{i(\delta,\varepsilon)}}(\tau) <
f_{\beta_{i(\delta,\varepsilon)}}(\tau) \\
  &< g_{i(\delta,\varepsilon)+1} < f_{\alpha_{i(\delta,\varepsilon)+1}}(\tau)
< f_\alpha(\tau)\} = {\frak a} \text{ mod } J
\endalign 
$$
\medskip

\noindent
hence ${\frak b}_\varepsilon \cap {\frak a}_{i(\delta,\varepsilon)} \notin
\emptyset$.  Moreover, ${\frak b}_\varepsilon \cap 
{\frak a}_{i(\delta,\varepsilon)}
\notin J$.  Now for each $\tau \in {\frak a}$ let $\varepsilon(\tau)$ be
$\sup\{\varepsilon < \kappa^+:\tau \in {\frak b}_\varepsilon \cap 
{\frak a}_{i(\delta,\varepsilon)}\}$ and let $\varepsilon(*) = 
\sup\{\varepsilon(\tau):\tau \in
{\frak a} \text{ and } \varepsilon(\tau) < \kappa^+\}$ so as $|{\frak a}|
\le \kappa$ clearly $\varepsilon(*) < \kappa^+$.  Let $\tau^* \in
{\frak b}_{\varepsilon(*)+1} \cap {\frak a}_{i(\delta,\varepsilon(*)+1)}$, so
$B = \{\varepsilon < \kappa^+:\tau^* \in {\frak b}_\varepsilon \cap
{\frak a}_{i(\delta,\varepsilon)}\}$ is unbounded in $\kappa^+,\langle
f_{\beta_{i(\delta,\varepsilon)}}(\tau^*):\varepsilon \in B \rangle$ is
strictly increasing (see clause $(c)$ above and the choice of
${\frak b}_\varepsilon$) and $\varepsilon \in B
\Rightarrow f_{\beta_{i(\delta,\varepsilon)}}(\tau^*) \in a_{\zeta(*)}$
(by the definition of ${\frak a}_{i(\delta,\varepsilon)}$, and $\zeta(*)$ as
$\zeta_{i(\delta,\varepsilon)} = \zeta(*)$).  We get contradiction to
$a \in \bold F(\kappa \cup \text{ Rang}(f_\alpha)) 
\Rightarrow |a| \le \kappa$.

So really we cannot carry the induction so we are stuck at some $i$.  If
$i=0$, or $i$ limit, or $i = j+1 \and \sup w_j < \chi$ we can find $g_i$
and then $\alpha_i$ and then $w_i$ as required.  So necessarily
$i =j+1,\sup(w_j) = \chi$.  Now if $\chi = \chi^*$, then this $w_j$ is as
required in the fact.  As pp$^+(\mu) = (\chi^*)^+$, the only case we cannot
have is when $\chi^*$ is singular.  Let $\chi^* = \underset{\varepsilon <
\text{ cf}(\chi^*)} {}\to \sup \chi_\varepsilon$ and 
$\chi_\varepsilon \in (\mu,\chi^*) \cap \text{ Reg}$ is 
(strictly) increasing with $\varepsilon$.  By \cite[Ch.II,\S1]{Sh:g} we 
can find, for each $\varepsilon < \text{ cf}(\chi^*),
{\frak a}_\varepsilon,J_\varepsilon,\bar f^\varepsilon = \langle
f^\varepsilon_\alpha:\alpha < \chi_\varepsilon \rangle$ as above, but in
addition
\medskip
\roster
\item "{$(*)$}"  $\bar f^\varepsilon$ is $\mu^+$-free i.e. for every $u \in
[\chi_\varepsilon]^\mu$ there is $\langle {\frak b}_\alpha:\alpha \in u
\rangle$ such that ${\frak b}_\alpha \in J_\varepsilon$ and for each $\tau
\in {\frak a}_\varepsilon,\langle f^\varepsilon_\alpha(\tau):\alpha
\text{ satisfies}:\tau \notin {\frak b}_\alpha \rangle$ is strictly 
increasing.
\endroster
\medskip

\noindent
So for every $a \in [\mu]^\kappa$ and $\varepsilon < \text{ cf}(\chi^*)$ 
we have

$$
\bigl\{\alpha < \chi_\varepsilon:\{\tau \in {\frak a}_\varepsilon:f_\alpha
(\tau) \in a\} \notin J_\varepsilon \bigr\} \text{ has cardinality } \le 
\kappa.
$$
\medskip

\noindent
Hence for each $a \in [\mu]^\kappa$

$$
\bigl\{(\varepsilon,\alpha):\varepsilon < \text{ cf}(\chi^*) \text{ and }
\alpha < \chi_\varepsilon \text{ and }\{\tau \in {\frak a}_\varepsilon:
f_\alpha(\tau) \in a\} \notin J_\varepsilon \bigr\}
$$
\medskip

\noindent
has cardinality $\le \kappa + \text{ cf}(\chi^*) = \text{ cf}(\chi^*)$
as for singular $\mu > \kappa \ge \text{ cf}(\mu)$ we have 
cf(pp$_\kappa(\mu) > \kappa$.
\medskip

\noindent
Define: $X = \{(\varepsilon,\alpha):\varepsilon < \text{ cf}(\chi^*),\alpha
< \chi_\varepsilon\}$

$$
\align
F'\bigl((\varepsilon,\alpha)\bigr) = \bigl\{ (\varepsilon',\alpha'):
&(\varepsilon',\alpha') \in X \backslash \{(\varepsilon,\alpha)\} \text{ and
for some} \\
  &d \in \bold F(\kappa \cup \text{ Rang}(f^\varepsilon_\alpha))
\text{ we have }
\{\tau \in {\frak a}_\varepsilon:f^{\varepsilon'}_{\alpha'}(\tau) \in d\} 
\notin J_{\varepsilon'} \bigr\}
\endalign
$$
\medskip

\noindent
so $F'\bigl((\varepsilon,\alpha)\bigr)$ is a subset of $X$ of cardinality
$< \text{ cf}(\chi^*)^+ + \theta < \chi^*$.

So by Hajnal's free subset theorem we finish (we could alternatively, for
$\chi^*$ singular, have imitated his proof). 
\hfill$\square_{\scite{2.14}}, \, \square_{\scite{2.15A}}$
\enddemo
\bigskip

\proclaim{\stag{2.16} Claim}  1) Assume
\medskip
\roster
\item "{$(a)$}"  $2^\lambda < 2^{\lambda^+} < 2^{\lambda^{+2}}$ and case B or
C of Fact \scite{2.10} for $\lambda$ occurs \nl
(so $\chi^*,T_\zeta$ are determined)
\sn
\item "{$(b)$}"  ${\frak K}$ is an abstract elementary class $LS({\frak K})
\le \lambda$ 
\sn
\item "{$(c)$}"  $K_{\lambda^{++}} \ne 0$,
\sn
\item "{$(d)$}"  ${\frak K}$ has amalgamation in $\lambda$ and in $\lambda^+$
\sn
\item "{$(e)$}"  in $K^3_\lambda$, the minimal triples are not dense. 
\ermn
\ub{Then}
\mr
\item "{$(\alpha)$}"  for each $\zeta < \chi^*$ for some
$M \in K_{\lambda^+}$ we have $|{\Cal S}_*(M)| \ge |\text{Lim}_{\lambda^+}
(T_\zeta)|$ \newline
(the tree from clause $(\zeta)$ of \scite{2.10}).
\ermn
2)  If $K$ satisfies (a)-(e) and is categorical in $\lambda^+$ or 
just has a universal member in $\lambda^+$, 
\ub{then} for some $M \in K_{\lambda^+}$ we have 
$|{\Cal S}_*(M)| = 2^{\lambda^+}$. \nl
3) If $(f)^+$, then $I(\lambda^{++},K) \ge 2^{\lambda^{++}}$ where
\mr
\item "{$(f)^+$}"  ${\frak K}$ is categorical in 
$\lambda$ and $\lambda^+$. 
\ermn
4)  Assume (a)-(e) of part (1) \ub{then} $\bold C^1_{{\frak K},\lambda}$ has
weaker $\lambda$-coding (if we have restricted to $(M,N,a)$ in
Definition of $\bold C$). 
\endproclaim
\bigskip

\remark{Remark}  Note that for \scite{2.16A} we do not use \scite{2.16}.
\endremark
\bigskip

\demo{Proof}  1) Let $\zeta < \chi^*$.  Let \newline
$T_\zeta = \dsize \bigcup_{\alpha < \lambda^+} T^\zeta_\alpha$ where 
$T^\zeta_\alpha$ are pairwise disjoint in $\alpha$, each $T^\zeta_\alpha$ has 
cardinality $\le \lambda$, \newline
$T^\zeta_0 = \{<>\}$ and $\eta \in T^\zeta_\alpha \and \beta <
\ell g(\eta) \Rightarrow \eta \restriction \beta \in \dsize \bigcup_{\gamma
< \alpha} T^\zeta_\gamma$, and \newline
$\eta \in T^\zeta_\alpha \Rightarrow \dsize \bigwedge_{\ell < 2} \eta 
\char 94 \langle \ell \rangle \in T^\zeta_{\zeta,\alpha +1}$  so 
$T_{\zeta,\alpha +1} = \{ \eta \char 94
\langle \ell \rangle:\eta \in T^\zeta_\alpha \text{ and } \ell < 2\}$.  For
$\eta \in T^\zeta_\delta,\delta$ a limit ordinal, necessarily both
$\ell g(\eta)$ and $\alpha(\eta) = \sup\{\gamma:\text{for some } 
\varepsilon < \ell g(\eta),\eta \restriction \varepsilon \in 
T^\zeta_\gamma\}$ are limit ordinals $\le \delta$.

Let $(M,N,a) \in K^3_\lambda$ be such that there is no minimal triple above
it. \newline
We now choose by induction on $\alpha < \lambda^+,\langle M^\zeta_\alpha,
M^\zeta_\eta,N^\zeta_\eta:\eta \in T^\zeta_\alpha \rangle$ such that:
\medskip
\roster
\item "{$(a)$}"  $(M^\zeta_\alpha,N^\zeta_\eta,a) \in K^3_\lambda$ and is
reduced if $\eta \in T^\zeta_\alpha,\alpha$ non-limit
\sn
\item "{$(b)$}"  $(M^\zeta_0,N^\zeta_{<>},a) = (M,N,a)$
\sn
\item "{$(c)$}"  $\nu \in T^\zeta_\beta,\eta \in T^\zeta_\alpha,\nu 
\triangleleft \eta,\beta < \alpha \Rightarrow 
(M^\zeta_\beta,N^\zeta_\nu,a) \le (M^\zeta_\alpha,N^\zeta_\eta,a)$ in the
order of $K^3_\lambda$
\sn
\item "{$(d)$}"  if $\delta$ is a limit ordinal then: $M^\zeta_\delta =
\dsize \bigcup_{\beta < \delta} M^\zeta_\beta$
\sn
\item "{$(e)$}"  if $\delta$ is a limit ordinal and $\eta \in T^\zeta_\delta$
then \newline
$N^\zeta_{\ell g(\eta)} = 
\dsize \bigcup_{\beta < \delta} N^\zeta_{\eta \restriction
\beta}$ hence $(M_{\alpha(\eta)},N^\zeta_\eta,a) \in K^3_\lambda$
\sn
\item "{$(f)$}"  if $\eta \in T^\zeta_\alpha$ then tp$(a,M_{\alpha +1},
N_{\eta \char 94 <0>}) \ne \text{ tp}(a,M_{\alpha +1},N_{\eta \char 94 <1>})$
\sn
\item "{$(g)$}"  $M^\zeta_\alpha \ne M^\zeta_{\alpha +1}$.
\endroster
\medskip

\noindent
There is no problem to carry the definition.  Let $M_\zeta =
\dsize \bigcup_{\alpha < \lambda^+} M^\zeta_\alpha \in K_{\lambda^+}$, and
for each $\nu \in \text{ Lim}_{\lambda^+}(T_\zeta)$ let
$N^\zeta_\nu = \dsize \bigcup_{\alpha < \lambda^+} N^\zeta_{\nu \restriction
\alpha}$, clearly $M_\zeta \le_{\frak K} N^\zeta_\nu$ and $a \in N^*_\nu$
and $\langle \text{tp}(a,M_\zeta,N^\zeta_\nu):\nu \in \text{ Lim}_{\lambda^+}
(T_\zeta) \rangle$ are pairwise distinct members of ${\Cal S}(M_\zeta)$.
This proves clause $(\alpha)$ of part (1).  \nl
2) This part follows by \scite{2.12}(2).  \nl
3) Now \scite{2.10} also applies to $(\lambda^+,\lambda^{++})$ in place of
$(\lambda,\lambda^+)$, so if for this case, clause $(A)_{\lambda^+}$ holds 
then by $(f)^+$ we get that \scite{2.14}(4) applies (with $(f)^+$ providing 
assumption (d)$^+$ there) hence we get
$I(\lambda^{++},K) = 2^{\lambda^{++}}$, so we can assume $(C)_{\lambda^+}
\vee (B)_{\lambda^+}$ occurs.

Now if WDmId$(\lambda^+)$ is not $\lambda^{++}$-saturated we get the desired
result as follows: by \scite{2.8} if $\neg(**)$ of \scite{2.7} holds and by 
\scite{2.6} we get a contradiction to $I(\lambda^+,K)=1$; if  
$(*)_\lambda$ (or $(*)'_\lambda$) of \scite{2.6} holds, but by \scite{2.7} 
one of the cases applies.
\newline
But as we are in case $(B)_\lambda$ or $(C)_\lambda$ (see (a) of 
\scite{2.16}(1)) by clause $(\iota)$ of \scite{2.10} we have 
$2^\lambda = \lambda^+,2^{\lambda^+} = \lambda^{++}$.  However, once we know 
$2^\lambda = \lambda^+$ we deduce that there is a model in $\lambda^+$ 
saturated over $\lambda$ so Claim \cite[3.16]{Sh:576} applies. \nl
4) Left to the reader as not used.   \hfill$\square_{\scite{2.16}}$
\enddemo
\bigskip

\proclaim{\stag{2.16A} Claim}  Assume
\medskip
\roster
\item "{$(a)$}"  $2^\lambda < 2^{\lambda^+} < 2^{\lambda^{++}}$
\sn
\item "{$(b)$}"  $K$ is categorical in $\lambda,\lambda^+$
\sn
\item "{$(c)$}"  $1 \le I(\lambda^{++},K) < 2^{\lambda^{++}}$.
\ermn
\ub{Then} above every $(M^*,N^*,a) \in K^3_\lambda$ there is a minimal triple.
\endproclaim
\bigskip

\demo{Proof}  Assume toward contradiction that above $(M^*,N^*,a) \in
K^3_\lambda$ there is no minimal type.  
If $2^\lambda = \lambda^+$, then there is a 
$M \in K_{\lambda^+}$ saturated over $\lambda$ hence we finish by 
\cite[3.16,At]{Sh:576}, possibility $(*)_2$.  So we can assume
$2^\lambda > \lambda^+$, hence \cite[6.3]{Sh:430} (with $\lambda^+$ here for
$\lambda$ there so $\mu$ there is $\le \lambda$ so $\delta < \lambda^+$ hence
$|T| \le |\delta| \le \lambda$ and let $\kappa = \text{ cf}
(\delta)$) there are $\kappa \le \lambda$ 
and \footnote{if $\kappa = \text{ Min}\{\sigma:2^\sigma > \lambda^+\}$, so
if $2^{< \sigma} \le \lambda$ then $({}^{\sigma >}2,\triangle)$ is okay,
otherwise ${}^{\sigma >}2 = \dbcu_{i < \lambda^+} T_i,|T_i| \le \lambda,T_i$
increasing with $i$ so for some $i,|\{\eta \in {}^\sigma 2:
\dsize \bigwedge_{\alpha < \sigma} \eta \restriction \alpha \in T_i\}| >
\lambda^+$} tree $T$ with
$\le \lambda$ nodes and $\kappa$ levels with $|\text{Lim}_\kappa(T)| >
\lambda^+$ hence for some $M \in K_\lambda,|{\Cal S}_*(M)| > \lambda^+$ (e.g.
by the proof of \scite{2.16}(1)).  If WDmId$(\lambda^+)$ is not
$\lambda^{++}$-saturated then in \scite{2.8} assumption (b) holds, and
assumptions (c) + (d) + (e) holds by the assumptions of the present claim
but not the conclusion, so (a) fails, that is $(**)_\lambda$ of \scite{2.7}
holds hence by \scite{2.7}, $(*)_\lambda$ of \scite{2.6} holds.  
But now \scite{2.6} contradicts clause (b) of the assumption, so we have to
assume that WDmId$(\lambda^+)$ is $\lambda^{++}$-saturated.  
Hence clause $(\iota)$ of \scite{2.10}, Case B does
not occur, hence Cases B,C of \scite{2.10} do not occur and hence Case A 
occurs.  So by \scite{2.14}(3) we get a contradiction to categoricity in
$\lambda^+$. \nl
${{}}$  \hfill$\square_{\scite{2.16A}}$
\enddemo
\bigskip

\centerline {$* \qquad * \qquad *$}
\bigskip

\proclaim{\stag{2.19} Claim}  We can prove \cite[4.2]{Sh:576} also for
$\lambda = \aleph_0$.
\endproclaim
\bigskip

\demo{Proof}  We ask: \newline
\underbar{Question 1}:  are there $M <_{\frak K} N$ in $K_\lambda$ such that
for no $a \in N \backslash M$ is tp$(a,M,N)$ minimal?

If the answer is yes, we can find $\langle M^1_i:i < \lambda^+ \rangle$
a representation of a model $M^1 \in K_{\lambda^+}$ such that:
$a \in M^1_{i+1} \backslash M^1_i \Rightarrow \text{ tp}(a,M^1_i,M^1_{i+1})$
is not minimal.  This implies $a \in M^1 \backslash M^1_i \Rightarrow
\text{ tp}(a,M^1_i,M^1)$ is not minimal (as for some $j \in [i,\lambda^+)$
we have $a \in M^1_{j+1} \backslash M^1_j$ so $(M^1_i,M^1_{j+1},a) \le
(M^1_j,M^1_{j+1},a)$ and the latter is not minimal).  But we can build another
representation $\langle M^2_i:i < \lambda^+ \rangle$ of $M^2 \in 
K_{\lambda^+}$ such that for each $i < \lambda^+$ for some 
$a \in M^1_{i+1} \backslash M^1_i$,
tp$(a,M^1_i,M^1_{i+1})$ is minimal (as there is a minimal triple).
So $M^1 \ncong M^2$.

So we assume the answer is no.
\medskip

\noindent
\underbar{Question 2}:  If $M \in K_\lambda,\Gamma \subseteq \Gamma^*_M =:
\{p \in {\Cal S} (M):p \text{ minimal}\}$ and $|\Gamma| \le \lambda$, is 
there $N$ such that: $M <_{\frak K} N \in K_\lambda$ and $N$ omit every 
$p \in \Gamma$?

If the answer to question 2 is yes, we can build $\langle M_\eta:\eta \in
{}^{\lambda^+}2 \rangle$ as in the proof of \scite{2.6}
(more exactly $\eta \triangleleft \nu \Rightarrow M_\eta \le_{\frak K} M_\nu,
M_\eta \in K_\lambda$) and we also have $\Gamma_\eta \subseteq \{p:\text{for
some } N \le_{\frak K} M_\eta,N \in K_\lambda \text{ and } p \in {\Cal S}(N)
\text{ is minimal not realized in } M_\eta\}$ have cardinality $\le \lambda,
\eta \triangleleft \nu \Rightarrow \Gamma_\eta \subseteq \Gamma_\eta$ and
there is $p \in \Gamma_{\eta \char 94 <1>}$ realized in
$M_{\eta \char 94 <0>}$ (and if you like also $p' \in \Gamma_{\eta \char 94
<0>}$ realized in $M_{\eta \char 94 <1>}$).  So by \cite[1.6]{Sh:576} we get
$I(\lambda^+,K) = 2^\lambda$.  So assume the answer is no
and for every $M \in K_\lambda$ let $\Gamma_M$ be a counterexample.  
Let $\langle M^1_i:i < \lambda^+ \rangle$, representing a
model $M^1 \in K_{\lambda^+}$ be such that $i < \lambda^+ \and p \in 
\Gamma_{M_i} \Rightarrow p$ realizes in $M$.  Now as in the proof of
saturated = model homogeneous (see \cite[0.21]{Sh:576}) we can prove $M^1$ is
saturated.  But this proves more than required: $|{\Cal S}(M^1_\ell)| \le 
\lambda^+$. \hfill$\square_{\scite{2.19}}$
\enddemo
\newpage

REFERENCES.  
\bibliographystyle{lit-plain}
\bibliography{lista,listb,listx,listf,liste}

\enddocument